\documentclass[11pt]{article}

\usepackage{amssymb,latexsym,amsmath,amsthm,graphicx,fullpage,titlesec}

\def\sqr#1#2{{\vcenter{\hrule height.#2pt              
     \hbox{\vrule width.#2pt height#1pt\kern#1pt
     \vrule width.#2pt}
     \hrule height.#2pt}}}
\def\square{\mathchoice\sqr{5.5}4\sqr{5.0}4\sqr{4.8}3\sqr{4.8}3}
\def\qed{\hskip4pt plus1fill\ $\square$\par\medbreak}

\def\bC{{\bf C}}
\def\bP{{\bf P}}

\def\bZ{{\bf Z}}

\def\cB{{\cal B}}

\def\cJ{{\cal J}}

\def\cR{{\cal R}}

\theoremstyle{plain}
\newtheorem*{theorem*}{Theorem}
\newtheorem{theorem}{Theorem}[section]
\newtheorem*{thmA}{Theorem A}
\newtheorem*{thm1}{Theorem 1}
\newtheorem*{thm2}{Theorem 2}
\newtheorem*{lemma*}{Lemma}
\newtheorem{lemma}[theorem]{Lemma}
\newtheorem*{proposition*}{Proposition}
\newtheorem{proposition}[theorem]{Proposition}
\newtheorem*{corollary*}{Corollary}
\newtheorem{corollary}[theorem]{Corollary}

\titleformat{\section}
  {\normalfont\bfseries}
  {\S\thesection.}{1em}{}

\begin{document}

\centerline{\bf Parabolic Bifurcations in Complex Dimension 2}
\bigskip



\centerline{Eric Bedford, John Smillie and Tetsuo Ueda}

\def\IMSmarkvadjust{0 pt}
\def\IMSmarkhadjust{0 pt}
\def\IMSmarkhpadding{0 pt}
\def\IMSpubltext{Published in modified form:}
\def\SBIMSMark#1#2#3{
 \font\SBF=cmss10 at 10 true pt
 \font\SBI=cmssi10 at 10 true pt
 \setbox0=\hbox{\SBF \hbox to \IMSmarkhpadding{\relax}
                Stony Brook IMS Preprint \##1}
 \setbox2=\hbox to \wd0{\hfil \SBI #2}
 \setbox4=\hbox to \wd0{\hfil \SBI #3}
 \setbox6=\hbox to \wd0{\hss
             \vbox{\hsize=\wd0 \parskip=0pt \baselineskip=10 true pt
                   \copy0 \break%
                   \copy2 \break%
                   \copy4 \break}}
 \dimen0=\ht6   \advance\dimen0 by \vsize \advance\dimen0 by 8 true pt
                \advance\dimen0 by -\pagetotal
	        \advance\dimen0 by \IMSmarkvadjust
 \dimen2=\hsize \advance\dimen2 by .25 true in
	        \advance\dimen2 by \IMSmarkhadjust

%
%
  \openin2=publishd.tex
  \ifeof2\setbox0=\hbox to 0pt{}
  \else 
     \setbox0=\hbox to 3.1 true in{
                \vbox to \ht6{\hsize=3 true in \parskip=0pt  \noindent  
                {\SBI \IMSpubltext}\hfil\break
                \input publishd.tex 
                \vfill}}
  \fi
  \closein2
  \ht0=0pt \dp0=0pt
 \ht6=0pt \dp6=0pt
 \setbox8=\vbox to \dimen0{\vfill \hbox to \dimen2{\copy0 \hss \copy6}}
 \ht8=0pt \dp8=0pt \wd8=0pt
 \copy8
 \message{*** Stony Brook IMS Preprint #1, #2. #3 ***}
}

\SBIMSMark{2012/1}{March 2012}{}

\bigskip

\setcounter{section}{-1}
\section{Introduction}   Parabolic bifurcations in one complex dimension demonstrate a wide variety of interesting dynamical phenomena [D, DSZ, L, Mc, S].  Consider for example the family of dynamical systems $f_\epsilon(z)=z+z^2+\epsilon^2$. When $\epsilon=0$ then $0$ is a parabolic fixed point for $f_0$. When $\epsilon\ne0$ the parabolic fixed point bifurcates into two fixed points. (The use of the term $\epsilon^2$ in the formula allows us to distinguish these two fixed points.)

We can ask how the dynamical behavior of $f_\epsilon$ varies with $\epsilon$. 
One way to capture this is to consider the behavior of dynamically significant sets such as the Julia set, $J$, and the filled Julia set, $K$, as functions of the parameter. 

\begin{theorem*}[{[D, L]}] 
The functions $\epsilon\mapsto J(f_\epsilon)$ and  $\epsilon\mapsto K(f_\epsilon)$ are discontinuous at $\epsilon=0$ when viewed as mappings to the space of compact subsets of $\bC$ with the Hausdorff topology. 
\end{theorem*}

In this paper we  consider parabolic bifurcations of families of diffeomorphisms in two complex dimensions. Specifically we consider a two variable family of diffeomorphisms $F_\epsilon: M\to M$ given locally by
$$F_\epsilon(x,y) = (x + x^2 + \epsilon^2+ \cdots, b_\epsilon y+\cdots)$$
where $|b_\epsilon|<1$, and the `$\cdots$' terms involve $x$, $y$ and $\epsilon$. When $\epsilon=0$ this map has the origin as a fixed point and the eigenvalues of the derivative at the origin are 1 and $b_0$. We say that $F_0$ is semi-parabolic  at the origin. 
 In [U1,2] it is shown that the set of points attracted to $O$ in forward time  can be written as $\cB\cup W^{ss}(O)$, where $\cB$ is a two complex dimensional basin of attraction and $W^{ss}(O)$ is the one complex dimensional strong stable manifold corresponding to the eigenvalue $b_0$.  The point $O$ is not contained in the interior of its attracting set, and we describe this by saying that the point is semi-attracting.   The set of points attracted to $O$ in backward time can be written as $\Sigma\cup O$ where $\Sigma$ is a one complex dimensional manifold.

A convenient two dimensional analog of the class of polynomial maps in one dimension is the family of polynomial diffeomorphisms of $\bC^2$. According to [FM] any dynamically interesting polynomial diffeomorphism is conjugate to a composition of generalized H\'enon maps; the degree 2 H\'enon map is given in (1.1). (For general discussions of such maps see [BS], [FS] and  [HO].)
Polynomial diffeomorphisms have constant Jacobian and to be consistent with the assumptions above we assume that the Jacobian is less than one in absolute value.  Analogs of the filled Julia set are the sets $K^+$ and $K^-$, consisting of points $p$ so that $F^n(p)$ remains bounded as $n\to\pm\infty$.  Analogs of the Julia set are the sets $J^\pm=\partial K^\pm$. We also consider $K=K^+\cap K^-$ and $J=J^+\cap J^-$. It is a basic fact that the one variable Julia set $J$ is the closure of the set of expanding periodic points.  We define $J^*$ to be the closure of the set of periodic saddle points. The set $J^*$ is contained in $J$ and has a number of other interesting characterizations: it is the Shilov boundary of $K$ and is the support of the unique measure of maximal entropy.  It is an interesting question whether these two sets are always equal.

\begin{thmA}
For $X=J^*$, $J$, $J^+$, $K$, $K^+$ the function $\epsilon\mapsto X(F_\epsilon)$ is discontinuous at $\epsilon=0$. For $X=J^-,K^-$ the function $\epsilon\mapsto X(F_\epsilon)$ is continuous. 
\end{thmA}

Our approach  follows the outlines of the approach of the corresponding result in one variable. In the one variable case the first step is to analyze certain sequences of maps $f_{\epsilon_j}^{n_j}$, where the parameter $\epsilon_j$ and the number of iterates $n_j$ are both allowed to vary.
The idea is the following. Let $p$ be a point in the basin of $0$ for $f_0$. When $\epsilon$ is small but non-zero the fixed point at $0$ breaks up into two fixed points.  As $n$ increases, the point $f^n_\epsilon(p)$ will come close to $0$ and may pass between these two fixed points and exit on the other side. Following standard terminology we refer to this behavior as ``passing through the eggbeater".  In this way it is possible to choose sequences $\epsilon_j$ and $n_j$ so that $f^{n_j}_{\epsilon_j}(p)$ will converge to some point on the other side of the eggbeater, in particular some point other than $0$.  The limit maps which arise this way have a convenient description in terms of Fatou coordinates of the map $f_0$. A {\sl Fatou coordinate} is a $\bC$-valued holomorphic map $\varphi$ defined on an attracting or repelling petal which satisfies the functional equation $\varphi(f(p))=\varphi(p)+1$.  There is an ``incoming'' Fatou coordinate $\varphi^\iota$ on the attracting petal and an ``outgoing'' Fatou coordinate $\varphi^o$ on the repelling petal.
Let $\tau_\alpha(\zeta)=\zeta+\alpha$ be the translation by $\alpha$, acting on $\bC$, and let $t_\alpha:=(\varphi^o)^{-1}\circ\tau_\alpha\circ\varphi^\iota$ for some $\alpha$.  Thus $t_\alpha$ maps the incoming petal to the repelling petal.

\begin{theorem*}[Lavaurs] 
If $\epsilon_j\to 0$ and $n_j\to\infty$ are sequences such that  $n_j-\pi/\epsilon_j\to\alpha$, then $\lim_{j\to\infty} f_{\epsilon_j}^{n_j}=t_\alpha$. 
\end{theorem*}

A sequence $(\epsilon_j,n_j)$ as in this Theorem will be called an $\alpha$-sequence.  Shishikura [S] gives a careful proof of this Theorem using the Uniformization Theorem in one dimension.  In Section 2, we re-prove the 1-dimensional result without using the Uniformization Theorem.  In Section 3 of this paper we prove the analogous result in two complex dimensions.

The existence of Fatou coordinates in the semi-parabolic case was established in [U1,2]. Let $\varphi^{\iota}:\cB\to\bC$ denote the Fatou coordinate on the attracting basin and $\varphi^{o}:\Sigma\to\bC$ the Fatou coordinate on the repelling leaf. Note that unlike the one variable case the function $\varphi^{\iota}$ has a two complex dimensional domain. In fact the map $\varphi^{\iota}$ is a submersion and defines a foliation whose leaves are described in Theorem 1.2. Define $T_\alpha:\cB\to \Sigma$ by the formula $T_\alpha=(\varphi^o)^{-1}\circ\tau_\alpha\circ\varphi^\iota$.  We introduce a useful normalization  (3.1), and  Theorem 3.1  shows that $F_\epsilon$ can be put in this form.  This simplifies the statement of the following:

\begin{theorem*}[3.9]
If $F_\epsilon$ satisfies (3.1), and if  $\epsilon_j\to 0$ and $n_j\to\infty$ are sequences such that $n_j-\pi/\epsilon_j\to\alpha$, then $\lim_{j\to\infty} F_{\epsilon_j}^{n_j}=T_{\alpha}$ at all points of ${\cal B}$. 
\end{theorem*}

We note that  we are taking very high iterates of a dissipative diffeomorphism, so the limiting map must have one-dimensional image. 


When $\epsilon$ is small, a point may pass through the eggbeater repeatedly.  We may use the map $T_\alpha$ to model this behavior.  In case $T_\alpha (p)$ happens to lie in ${\cal B}$, we may define the iterate $T_\alpha^2(p)$.    A point for which $T^n_\alpha$ can be defined for $n$ iterations corresponds to a point which passes through the eggbeater $n$ times.

Following the approach of [D, L] in one dimension we may introduce sets $J^*(F_0, T_\alpha)$ and $K^+(F_0, T_\alpha)$ which reflect the behavior points in $\cB$ under the maps $F_0$ and $T_\alpha$.

\begin{thm1}
Suppose that $F_\epsilon$ is normalized as in (3.1).   
If  $\epsilon_j$ is an $\alpha$-sequence, then
$$\liminf_{j\to\infty} J^*(F_{\epsilon_j})\supset J^*(F_0,T_\alpha),$$
where `lim inf' is interpreted in the sense of Hausdorff  topology.
\end{thm1}

Though we have stated Theorem 1 for polynomial diffeomorphisms in fact the definition of the set $J^*(F_\epsilon)$ makes sense for a general holomorphic diffeomorphism and Theorem 1 is true in this broader setting.

\bigskip
\centerline{\includegraphics[height=2.3in]{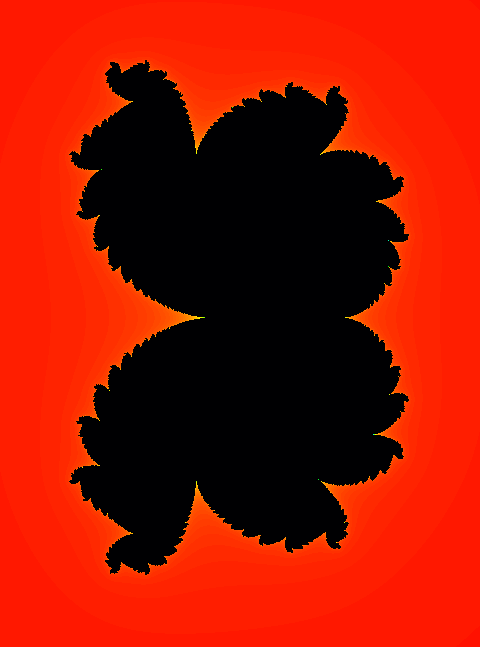}   \ \ \  \ \ \   \includegraphics[height=2.3in]{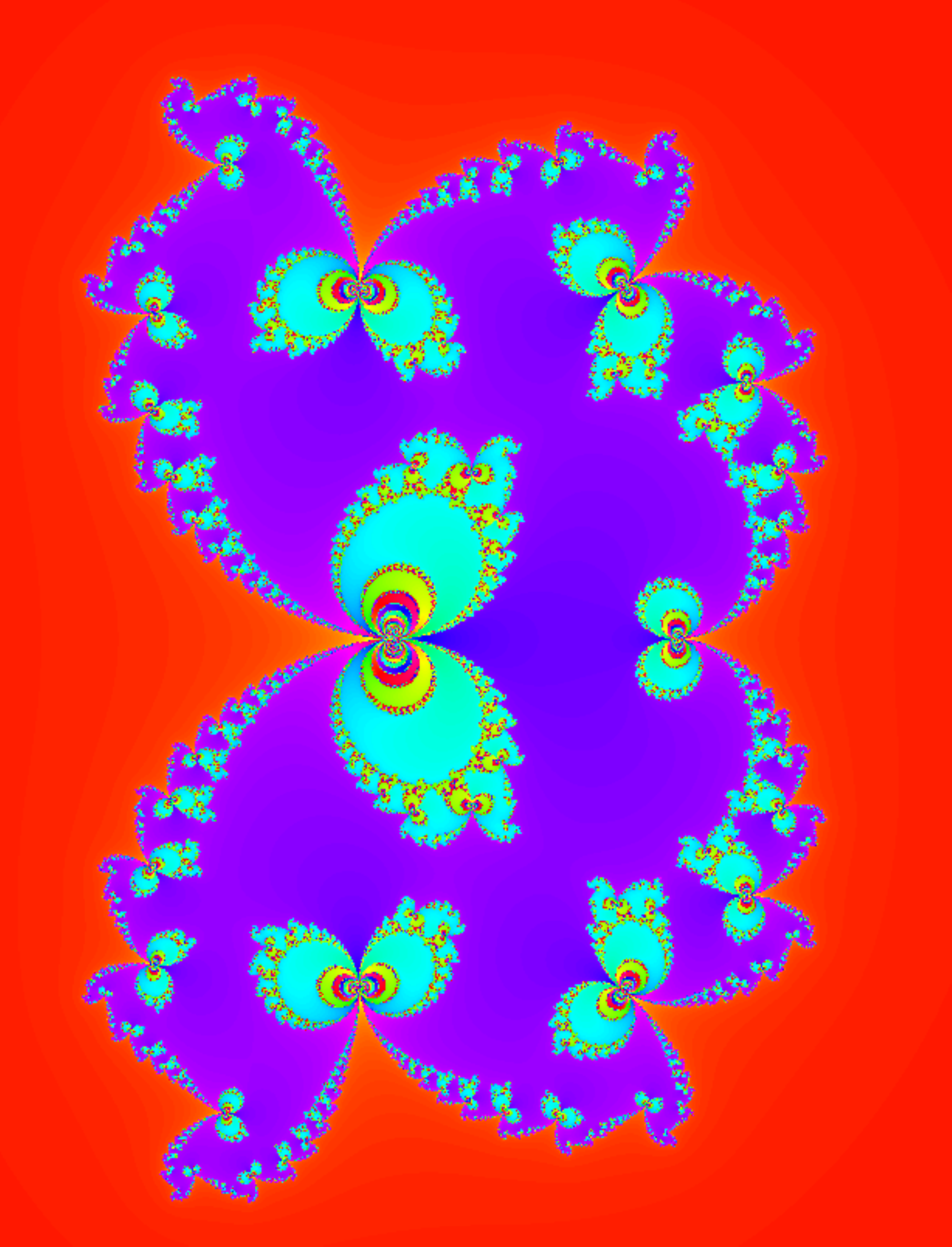} }

\centerline{Figure 1.  Discontinuity of $\epsilon\mapsto K^+(F_{a,\epsilon})$ illustrated  by complex linear slices in ${\bf C}^2$.}  

\centerline{\ \ \ \ \ \ \ \ \ $F_{a,\epsilon}$ is given by equation (1.1) with  $a=.3$; $\epsilon=0$ (left), $a=.3$, $\epsilon=.05$ (right).}  

\bigskip

\begin{thm2}
Suppose that $F_\epsilon$ is conjugate to a composition of generalized H\'enon maps, and $F_\epsilon$ is normalized as in (3.1).  If  $\epsilon_j$ is an $\alpha$-sequence, then we have
$$\cB\cap \limsup_{j\to\infty}K^+(F_{\epsilon_j})\subset K^+(F_0,T_\alpha).$$ 
\end{thm2}

If the function $\epsilon\mapsto J^*(F_\epsilon)$ were continuous at $\epsilon=0$, then the limit of $J^*(F_{\epsilon_j})$ along an $\alpha$-sequence would be independent of $\alpha$ and would be equal to $J^*(F_0)$. Theorem 1 implies that $J^*(F_0)$ would have to contain every set $\cJ_\alpha^*$. Theorem 2 implies that $J^*(F_0)$ would have to be contained in every set 
$K^+(F_0,T_\alpha)$. Our next result shows that these conditions are incompatible. 

\begin{theorem*}[4.4] 
For each $p\in\cB$ there are constants $\alpha$ and  $\alpha'$ such that $p\in J^*(F_0,T_\alpha)$ but $p\notin K^+(F_0,T_{\alpha'})$. 
\end{theorem*}

We can use Theorem 4.4 to prove the discontinuity statement of the maps $\epsilon\mapsto J^*(F_\epsilon)$ and $\epsilon\mapsto K^+(F_\epsilon)$, but in fact the same argument shows the discontinuity of any dynamically defined set $X$ which is sandwiched between $J^*$ and $K^+$. Using this idea, we now give a proof of Theorem~A. 
\medskip

\noindent{\bf Proof of Theorem A.} We begin by proving the statement concerning discontinuity.
Let $X$ be one of the sets $J$, $J^*$, $J^+$, $K$, $K^+$. Assume that the function $\epsilon\mapsto X(F_\epsilon)$ is continuous at $\epsilon=0$. Choose a $p$ in $\cB$ and let $\alpha$ and $\alpha'$ be as in Theorem 3. Let $\epsilon_j$ be an $\alpha$ sequence and let $\epsilon'_j$ be an $\alpha'$ sequence. Since $J^*(F_\epsilon)\subset X(F_\epsilon) \subset K^+(F_\epsilon)$ we have that $J^*(F_0,T_\alpha)\subset \liminf_{j\to\infty} J^*(F_{\epsilon_j})\subset \liminf_{j\to\infty} X(F_{\epsilon_j})=X(F_0)$ by Theorem 1 and $K^+(F_0,T_{\alpha'})\supset \liminf_{j\to\infty} K^+(F_{{\epsilon'}_j})\supset \liminf_{j\to\infty} X(F_{{\epsilon'}_j})=X(F_0)$ by Theorem 2. This gives $J^*(F_0,T_\alpha)\subset K^+(F_0,T_\alpha')$. On the other hand $p\in J^*(F_0,T_\alpha)$ but $p\notin K^+(F_0,T_\alpha')$ by Theorem 3 so we arrive at a contradiction.

The fact that $J^-$ and $K^-$ vary continuously follows from the fact that for polynomial diffeomorphisms which contract area the sets $J^-$ and $K^-$ are equal (see [FM]). We combine this with the facts from Proposition 4.7  that $J^-$ varies lower semi-continuously and $K^-$ varies upper semi-continuously.

\medskip

Computer pictures illustrate the behaviors
described in Theorems~1 and 2.
Consider the family of quadratic H\'enon diffeomorphisms of ${\bf C}^2$:
$$ F_{a,\epsilon}(x,y) =((1+a)x-ay+x^2+\epsilon^2, x+\epsilon^2).  \eqno{(1.1)} $$
The parameters are chosen so that $F$ is semi-parabolic
when $\epsilon=0$; the origin $O$ is the unique fixed point and has multipliers 1 and $a$.  Figure 1 shows the slice $K^+\cap T$, where $T$ is the complex line passing through $O$ and corresponding to the eigenvalue 1 when $\epsilon=0$.   We color points according to the value of the Green function.  The set $K^+=\{G^+=0\}$ is colored black.  It is hard to see black in the right hand of Figure 1 because the set $T\cap K^+$ is small.  But we note that $G^+$ is harmonic where it is nonzero, so points of $T\cap K^+$ must be present in the apparent boundaries between regions of different color.
In the perturbation shown in Figure 1, there is not much change to the ``outside'' of $K^+$, whereas the ``inside'' shows the effect of an ``implosion.''  Further discussion of the figures in this paper is given at the end of \S1.

\bigskip\noindent{\bf Acknowledgements.}  We thank M.~Shishikura for his generous advice and encouragement throughout this work.  This project began when E.B. and J.S. visited RIMS for a semester, and they wish to thank Kyoto University for its continued hospitality.  We also wish to thank H.~Inou and X.~Buff.

\section{Fatou coordinates and transition functions.}  
Let $M$ be a 2-dimensional complex manifold, and let $F$ be an automorphism of $M$.  Let $O$ be a fixed point which  is a semi-attracting, semi-parabolic.  By [U1] we may choose $i$ and $j$ and change coordinates so that $O=(0,0)$, and $F$ has the form
\begin{equation}
\begin{split}
x_1 & = x + a_2 x^2 + \dots + a_i x^i + a_{i+1}(y) x^{i+1} + \dots \\
y_1 &= by + b_1 xy+\dots + b_j x^jy+ b_{j+1}(y)x^{j+1}+\dots
\end{split}
\tag{1.2}
\end{equation}

We will suppose that $a_2\ne0$, and thus by scaling coordinates, we may assume $a_2=1$.  (In the case where $a_2=\dots=a_m=0$, $a_{m+1}\ne0$, the results analogous to [U1] have been treated by Hakim [H].)  We may choose coordinates so that the local stable manifold is given by $W^s_{loc}(O)=\{x=0,|y|<1\}$.   For $r,\eta_0>0$, we set $B^\iota_{r,\eta_0}=\{|x+r|<r,|y|<\eta_0\}$.  If we take $r,\eta_0$ small, then the iterates $F^n\overline B^\iota_{r,\eta_0}$ shrink to $O$ as $n\to\infty$.  Further, $B^\iota_{r,\eta_0}$ plays the role of the ``incoming petal'' and is a base of convergence in the sense of [U1], which is to say that $\cB:=\bigcup_{n\ge0}F^{-n}B_{r,\eta_0}^\iota$ is the set of points where the forward iterates converge locally uniformly to $O$.




With $a_3$ as in (1.2), we set $q = a_3-1$ and choosing the principal logarithm, we set
$$w^\iota(x,y):= -{1\over x} - q \log(-x). \eqno(1.3)$$
It follows (see [U1]) that for $p\in\cB$ the limit 
$$\varphi^\iota(p) = \lim_{n\to\infty}(w^\iota(F^n(p)) - n)$$
converges to an analytic function $\varphi^\iota:\cB\to\bC$ satisfying the property of an Abel-Fatou coordinate:  $\varphi^\iota\circ F = \varphi^\iota +1$.  Further, 
$$\varphi^\iota (x,y)-w^\iota(x)=B(x,y) \eqno(1.4)$$ 
where $B$ is bounded on $B_{r,\eta_0}^\iota$, and in fact this condition defines $\varphi^\iota$  up to additive constant. 

%


Let us note a  result from [U1]:
\begin{theorem}
There is an entire function $\Phi_2(x,y)$ such that $\Phi=(\varphi^\iota,\Phi_2)$ is biholomorphic 
$\Phi=(X,Y):\cB\to\bC^2$, and 
$f$ corresponds to $(X,Y)\mapsto(X+1,Y)$. 
\end{theorem}

By Theorem 1.1,  $\varphi^\iota$ has nonvanishing gradient and thus defines a foliation of $\cB$ whose leaves are closed complex submanifolds which are holomorphically equivalent to $\bC$.  We conclude from (1.4) that  if we take $r$ small, then for fixed $|y_0|<\eta_0$,  $x\mapsto \varphi^\iota(x,y_0)$ is univalent on $\{|x+r|<r\}$, and the image is approximately $G_r:=\{w^\iota(|x+r|<r)$.  It follows that there is a domain $G\subset{\bf C}$ such that if $\xi\in G$, then there is an analytic function $\psi_\xi(y)$ for $|y|<\eta_0$ such  that
$$\{(x,y)\in B_{r,\eta_0}:\varphi^\iota(x,y)=\xi\} = \{x=\psi_\xi(y):|y|<\eta_0\}.$$
We may choose $0<r_1<r_2$ such that $B_{r_1,\eta_0}$ is contained in the union of such graphs, and each of these graphs is contained in $B_{r_2,\eta_0}$.  We use this to prove that the fibers $\{\varphi^\iota=const\}$ are strong stable manifolds in the sense of exponential convergence in (1.5), whereas the convergence in the parabolic direction is quadratic (cf.\ [Mi, Lemma 10.1]).  

\medskip
\vbox{\centerline{ \includegraphics[height=2.4in]{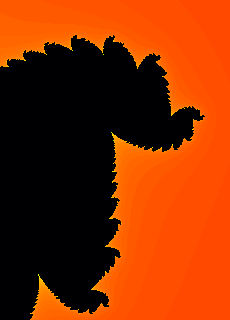}  \ \ \ \ \ \  \includegraphics[height=2.4in]{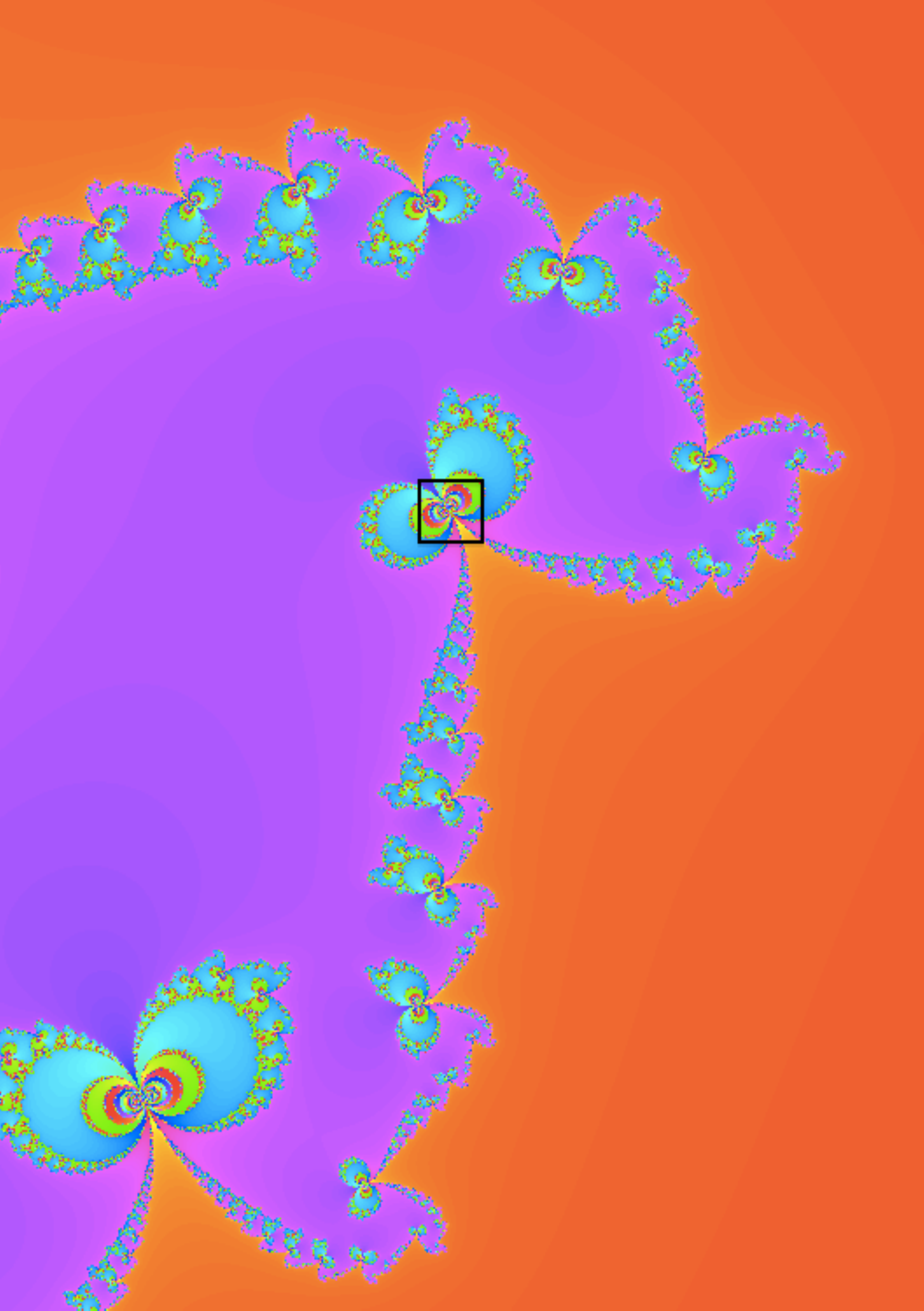}   } 

\centerline{Figure 2a.  Unstable slices of $K^+$ for $F_{a,\epsilon}$. $a=.3$, $\epsilon=0$ (left); $a=.3$, $\epsilon=.05$ (right).}}

\bigskip

\bigskip
\centerline{  \includegraphics[height=2.7in]{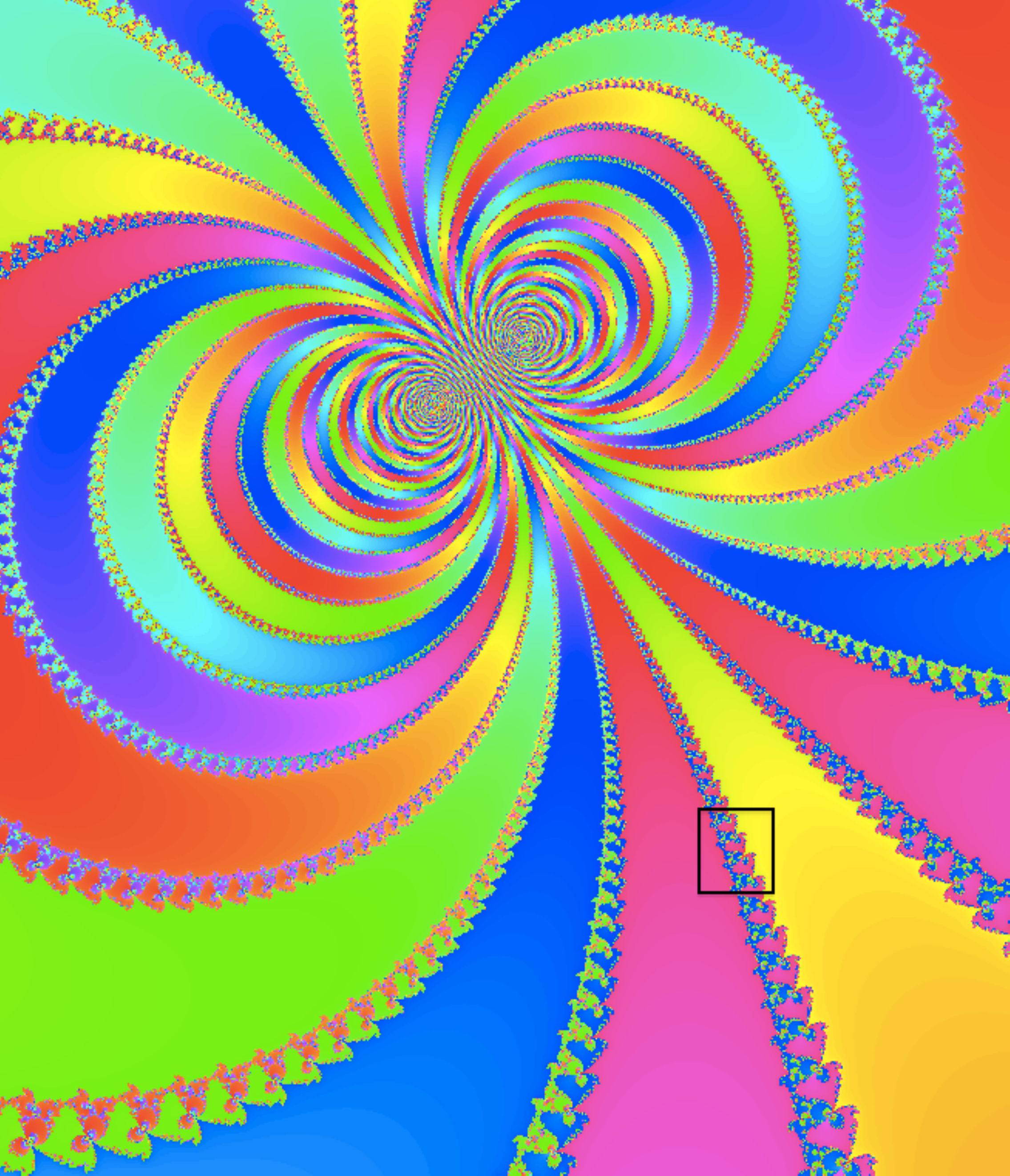}  \ \ \ \ \ \  \includegraphics[height=2.7in]{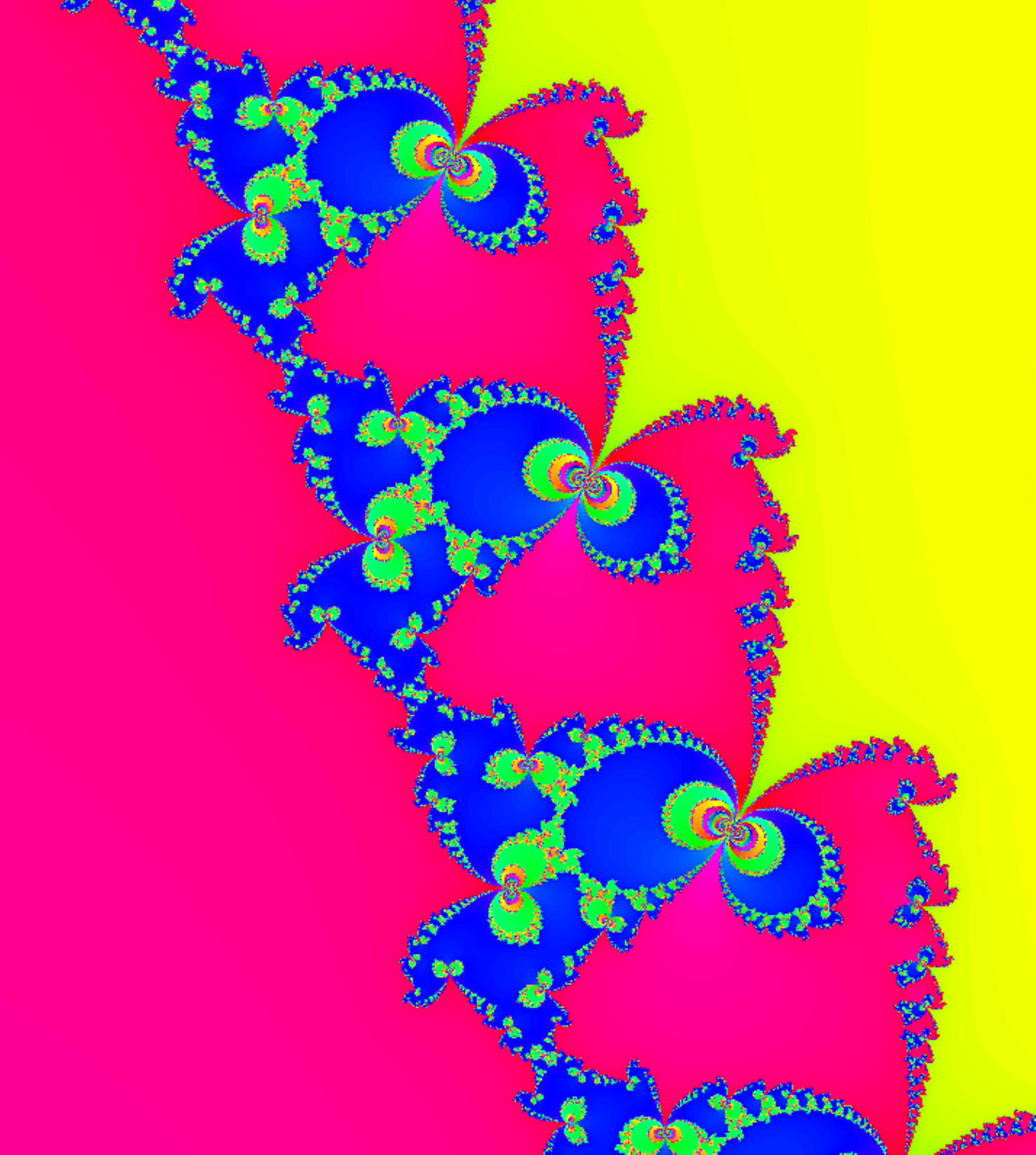} }

\centerline{Figure 2b.  Zooms of Figure 2a.}

\medskip

\begin{theorem}
For $p_1,p_2\in{\cal B}$ such that $p_1\ne p_2$ and $\varphi^\iota(p_1)=\varphi^\iota(p_2)$, we have
$$\lim_{n\to+\infty}{1\over n}\log\,  {\rm dist}(F^np_1,F^np_2)=\log|b|.\eqno(1.5)$$
Conversely, if $\varphi^\iota(p_1)\ne\varphi^\iota(p_2)$, then $\lim_{n\to+\infty}\left( n^2\cdot {\rm dist}(F^np_1,F^np_2)\right)\ne0$.
\end{theorem}

\noindent{\it Proof.  }  If we iterate the points forward, they will enter the set $B_{r,\eta_0}$, so we may assume that $p_1,p_2\in B_{r,\eta_0}$.  If $\varphi^\iota(p_1)=\varphi^\iota(p_2)$ we may assume that they are contained in a graph $\{x=\psi_\xi(y):|y|<\eta_0\}$.  The behavior in the $y$-direction is essentially a contraction by a factor of approximately $|b|$, so the distance from $F^np_1$ to $F^np_2$ is essentially contracted by $|b|$, which gives the first assertion.  

For the second assertion, we note that along a forward orbit we have  
$$\varphi^\iota(x,y) = - {1\over x} -q\log(-x) + \gamma+ o(1),$$ 
where $o(1)$ refers to a term which vanishes as the orbit tends to $O$.   Without loss of generality, we may suppose that $\gamma=0$.  From this we find  
$$x -q x^2  \log(-x) + \cdots = - {1\over \varphi^\iota}.$$
Now we substitute this expression into itself and obtain
$$x =  - {1 \over \varphi^\iota} + {q\over (\varphi^\iota)^2}\log\left({1\over\varphi^\iota}\right) + \cdots.$$ 
If we write $\varphi^\iota(x_i,y_i) = c_i$ for $i=1,2$, then $\varphi^\iota(f^n(x_i,y_i))=c_i+n$, and $f^n(x_i,y_i) = (x_{i,n},y_{i,n})$ satisfies 
$$x_{i,n}= - {1\over c_i+n} + {q\over (c_i+n)^2}\log(c_i+n)^{-1} + \cdots$$
Thus $x_{1,n}-x_{2,n}= (c_1-c_2)/\left( (c_1+n)(c_2+n)\right) + O(n^{-3}\log n)$, so $\lim_{n\to\infty} n^2|x_{1,n}-x_{2,n}| = |c_1-c_2|$. \qed

We may also define the asymptotic curve
$$\Sigma:=\{p\in M-\{O\}: f^{-n}p\to O\ {\rm as\ }n\to\infty)\}.  \eqno(1.6)$$
This is a Riemann surface which is equivalent to $\bC$.  Let us define $B^o_r:=\{|x-r|<r,|y|<\eta_0\}$, which is the analogue of the ``outgoing petal.''   By [U2], there is a component $\Sigma_0$ of $\Sigma\cap B^o_r$ such that $\overline\Sigma_0$ is a smooth graph $\{\psi(x)=y, |x-r|<r\}$, and $\psi(0)=0$.  Thus $O$ is in the boundary of a smooth piece of $\Sigma$.  If $\psi$ extended analytically past $x=0$, then $\Sigma\cong\bC$ would be contained in a larger complex manifold, which would have to be $\bP^1$.  Thus $M$ would  contain a compact, complex curve.  Stein manifolds (${\bf C}^2$, for instance) do not have such curves, so we have:

\begin{proposition}
If $M$ is Stein, then $\Sigma$ cannot be extended analytically past $O$.
\end{proposition}


\medskip
\noindent{\bf Examples.} The first example is the product $M_0=\bP^1\times \bC$.  Let $F$ act as translation on $\bP^1\times\{0\}$ with fixed point $O = (\infty,0)\in\bP^1\times\{0\}$, and let $F$ multiply the factor $\bC$ by $b$.  Then $\Sigma=\bP^1\times \{0\} - O$, and $\cB=\Sigma\times\bC$, so we see that $\Sigma\subset\cB$.

For the second example, we start with the linear map  $L(x,y) = (b(x + y), by)$ on $\bC^2$, so the $x$-axis $X=\{y=0\}$ is invariant.  The origin is an attracting fixed point, and we let  $M_1$ denote $\bC^2$ blown up at the origin.  Thus $L$ lifts to a biholomorphic map of $M_1$.  We write the exceptional fiber as $E$ and note that  $E$ is equivalent to $\bP^1$, and $L$ is equivalent to translation on $E$.  The fixed point of $L|_E$ is $E\cap X$.  We have  $\Sigma=E-X$ and $\cB=M_1-X$, and $\cB$ contains $\Sigma$.  The second example is different from the first because $E$ has negative self-intersection.  

Both $M_0$ and $M_1$ fail to be Stein because they contain compact holomorphic curves.  Similar examples can be constructed for all of the Hirzebruch surfaces.
\bigskip

\centerline{\includegraphics[width=3.2in]{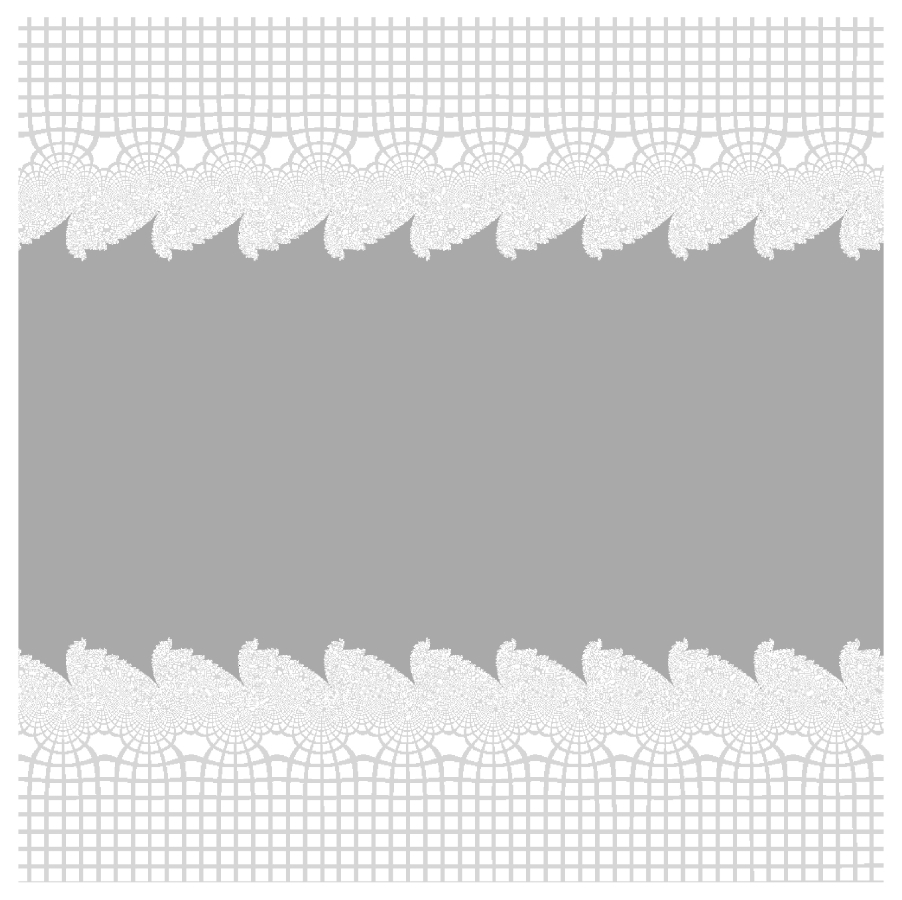}   \includegraphics[width=3.2in]{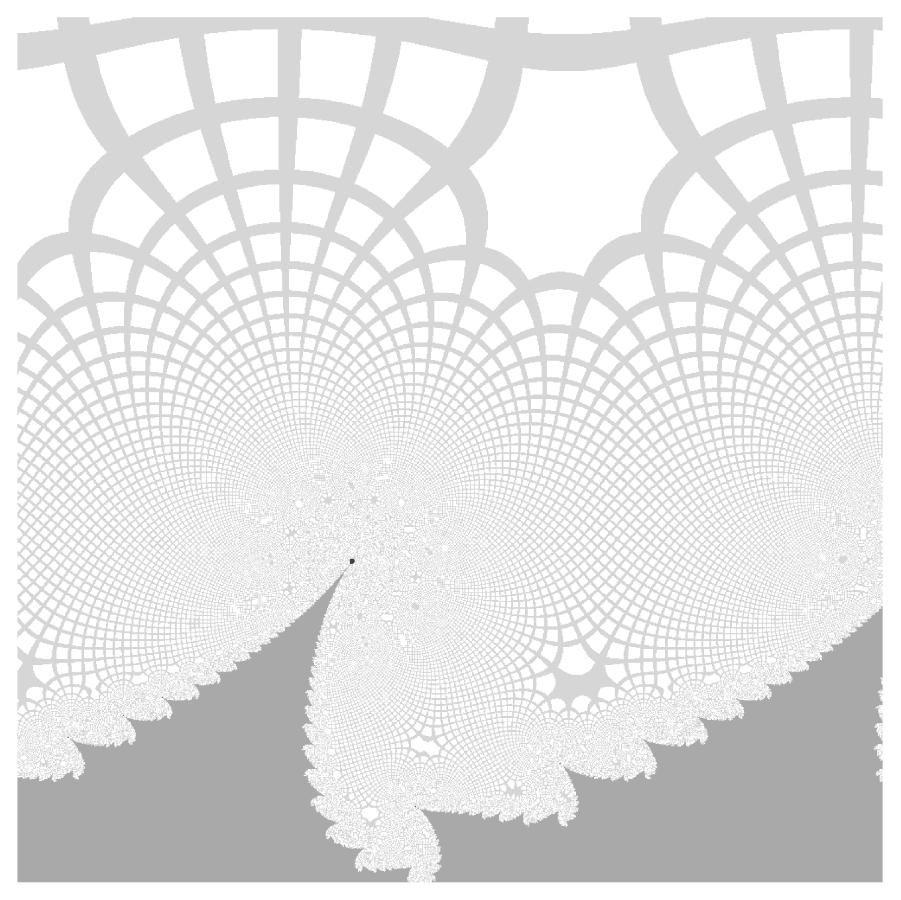}  }

\centerline{Figure 3.  Basin ${\cal B}$ in Fatou coordinates; $a=.3$; 10 periods (left) and detail (right).} 
\medskip

\medskip
In order to define the outgoing Fatou coordinate, we set
$$w^o(x,y):= -{1\over x} - q \log(x)$$
and we define a map $\varphi^o:\Sigma\to\bC$ by setting
$$\varphi^o(p) = \lim_{n\to\infty}(w^o(F^{-n}(p)) + n)$$
This satisfies $\varphi^o(f)=\varphi^o+1$, and it is a bijection (see [U2]). 
For $\alpha\in\bC$ we define the translation $\tau_\alpha:\bC\to\bC$ by $\tau_\alpha(\zeta)=\zeta+\alpha$.  We define
$$T_\alpha:= (\varphi^o)^{-1}\circ \tau_\alpha\circ\varphi^\iota:\cB\to\Sigma$$
We have
$$F\circ T_\alpha=T_{\alpha+1}=T_\alpha\circ F.$$
Since $\varphi^\iota$ and $\varphi^o$ are defined up to additive constants, the family $\{T_\alpha\}$ is independent of the choice of $\varphi^\iota$ and $\varphi^o$.

Let us set $\Omega:=\varphi^o(\cB\cap\Sigma)\subset\bC$.  Thus we have a family of maps $h_\alpha:\Omega\to\bC$ given by
$$h_\alpha:=\tau_\alpha\circ \varphi^\iota\circ(\varphi^o)^{-1}$$
These maps  all agree up to an additive constant, so the  maps  $\{h_\alpha\}$ all have the same set of  critical points; changing $\alpha$ serves to change the critical values.  The critical points correspond to the points $\zeta_c$ where $(\varphi^o)^{-1}(\zeta_c)$ is a point of tangency between $\Sigma$ and the strong stable fibration $\{\varphi^\iota=const\}$.

Since we may iterate the map $h_\alpha$ as long as the point stays inside $\Omega$, this yields a family of partially defined dynamical systems.  Each map $h_\alpha$ satisfies $h_\alpha(\zeta+1)=h_\alpha(\zeta)+1$.  For  $R>0$, let us write $\Omega^\pm_R:=\{\zeta\in\bC: \pm\Im\zeta>R\}$, and choose $R$ large enough that $\Omega^\pm_R\subset\Omega$.  On $\Omega^\pm_R$ we have
$$h_{\alpha}(\zeta) = \zeta + \alpha +  c_0^\pm + \sum_{n>0}c^\pm_{ n}e^{\pm 2n\pi i\zeta}.$$ 
In particular $h_\alpha$ is injective on $\Omega^\pm_R$ if $R$ is sufficiently large.  Since $h_\alpha$ is periodic, it defines a map of the cylinder ${\bf C}/{\bf Z}$; we see that the upper (resp.\ lower) end of the cylinder will be attracting if $\Im(\alpha + c_0^+)>0$ (resp.\ $\Im(\alpha+c_0^-)<0$).

In the construction of the Fatou coordinates we have
\begin{align*}
&\varphi^\iota(x,y)=-{1\over x} - q\log(-x) +o(1), \ \ (x,y)\in{\cal B}_{\rm loc},\cr
& \varphi^o(x,y)=-{1\over x} - q\log(x) + o(1), \ \ (x,y)\in\Sigma_{\rm loc}.
\end{align*}
Thus if we compare the values of $\log$ at the upper and lower ends of the cylinders and use this in the formula for $h_\alpha$, we find $c_0^\pm=\pm \pi i q$, which gives the normalization $c_0^++c_0^-=0.$   For comparison, we note that Shishikura [S] uses the normalization  $c_0^+=0$.   In the case of the semi-parabolic map (1.1) with $\epsilon=0$, we find that if we normalize the form (1.2) so that $a_2=1$, then we have  $a_3=2a/(a-1)$ and $q=a_3-1$, so
$$c_0^+ = \pi i{a+1\over a-1}  \eqno(1.7)$$

%



We make an elementary observation:
\begin{proposition}
If $M$ is Stein, then $\cB$ is a component of normality of the family $\{f^n:n\ge0\}$, and thus $\cB$ is polynomially convex in $M$.  Since $\Sigma$ is simply connected, it follows that every component of $\cB\cap\Sigma$ is simply connected.
\end{proposition}

Let $\Omega^\pm$ denote the component of $\Omega$ which contains $\Omega^\pm_R$.

\begin{theorem}
Assume that  $M$ is Stein, and $\Omega\ne{\bf C}$.  If $\Omega'$ is a connected component of $\Omega$, then the function $h_\alpha|_{\Omega'}$  cannot be continued analytically over any boundary point of $\Omega'$.  In particular, since $\Omega^\pm\ne{\bf C}$, the derivative $h_\alpha'$ is nonconstant on both $\Omega^+$ and $\Omega^-$, and there exist points in both of these sets where $| h'_\alpha|<1$ and where $|h'_\alpha|>1$.
\end{theorem}

\noindent{\it Proof.}  Let us fix a boundary point $\zeta_0\in\partial\Omega'$.  Let $\Delta$ be a disk containing $\zeta_0$, and let $\Delta_1$ be a component of $\Delta\cap\Omega'$.  We will show that $h_\alpha$ is not bounded on $\Delta_1$.  We know that $(\varphi^o)^{-1}:{\bf C}\to\Sigma$ is entire, so  $(\varphi^o)^{-1}(\zeta)\to \Sigma\cap\partial{\cal B}$ as $\Omega\ni\zeta\to\partial\Omega$.  If  $\Phi=(\varphi^\iota,\Phi_2)$ is the map from Theorem 1.1, then we have $||\Phi((\varphi^o)^{-1}(\zeta))||^2=|\varphi^\iota((\varphi^o)^{-1}(\zeta))|^2+|\Phi_2((\varphi^o)^{-1}(\zeta))|^2\to \infty$ as $\zeta\in\Delta_1$ approaches $\partial\Omega\cap\overline{\Delta_1}$.  If $h_\alpha$ is bounded near $\zeta_0$, then so is $\varphi^\iota$.  It follows that $|\Phi_2(\zeta)|\to\infty$ as $\zeta\to\partial\Delta_1$ near $\zeta_0$.  But this is not possible, since by Proposition 1.4, $\Omega'$ is simply connected, so $\partial\Delta_1$ has no isolated boundary points.  

Thus $h_\alpha$ cannot be constant on any component of $\Omega$.  In particular, the derivative is not constantly zero on $\Omega^\pm$.  Since $\lim_{\zeta\to\infty}h'_\alpha=1$, by the Maximum Principle there must be points near infinity where $|h'_\alpha|>1$ and where $|h'_\alpha|<1$.  \qed


We let $\bar\Omega$ denote the image of $\Omega$ in the cylinder $\bC/\bZ$.  Then $h_\alpha$ passes to an analytic map $\bar h_\alpha:\bar\Omega\to\bC/\bZ$.  Further, $\bar h_\alpha$ extends analytically past each of the ends of the cylinder.  For instance, at the upper end of the cylinder, $\bar h_\alpha$ is analytic, as a function of the variable $z = e^{2 \pi i\zeta}$, in a neighborhood of $z=0$.

Let us close this section with the comment that certain aspects of this construction are local at $O$.  In case $F$ is defined in a neighborhood $U$ of $O$, we may define the local basin 
$$\cB_{loc}:=\{p: f^np\in U\ \forall n\ge0, f^np\to O{\rm\ locally\ uniformly\ as\ } n\to\infty\},$$ 
as well as the local asymptotic curve $\Sigma_{loc}$.  Similarly, we have Fatou coordinates $\varphi^{\iota}$ and $\varphi^o$ on $\cB_{\rm loc}$ and $\Sigma_{\rm loc}$.    In this case there is an $R$ such that 
$$\varphi^\iota(\cB_{\rm loc})\supset\{\zeta\in\bC: -\Re\zeta +R<|\Im\zeta|\}, \ \ \   \varphi^o(\Sigma_{\rm loc})\supset\{\zeta\in\bC: \Re\zeta +R<|\Im\zeta|\}.$$
We define $W_R:= \{\zeta\in\bC:|\Re\zeta| + R<|\Im\zeta|\}$, so for $R$ sufficiently large,  
$$\varphi^{\iota/o}(\cB_{\rm loc}\cap \Sigma_{\rm loc})\supset W_R,$$
and possibly choosing $R$ even larger, $h_\alpha=\tau_\alpha\circ \varphi^\iota\circ H$ is defined as a map of  $W_R$ to $\bC$.  Note that we have $h_\alpha(\zeta+1)= h_\alpha(\zeta)+1$ for $\zeta\in W_R$ such that both sides of the equation are defined.   If we shrink the domain $U$ of $F$, we may need to increase $R$, but the germ of $h_\alpha$ at infinity is unchanged.  Let  $h_\alpha^\bullet$ we denote the germ at infinity of $h_\alpha$ on $W_R$.  It is evident that:

\begin{theorem}
If $F$ and $F'$ are locally holomorphically conjugate at $O$, then there is a translation on $\bC$ which conjugates the families of germs $\{h_\alpha^\bullet\}$ to $\{{h'}_\alpha^\bullet\}$.
\end{theorem}

\bigskip\noindent{\bf Graphical representation.}   In Figure 1, we saw slices of sets by planes.  For an invariant picture, we may slice  by
unstable manifolds of periodic saddle points.
If $Q$ is a periodic point of saddle type, then the unstable manifold $W^u(Q)$ may be
uniformized by $\bC$ so that $Q$ corresponds to $0\in \bC$.
The restriction of $F$ to $W^u(Q)$   corresponds to a linear map of $\bC$ 
in the uniformizing coordinate, so the slice picture is self-similar.
The unstable slice picture cannot be taken at the fixed point $O$ 
when $\epsilon=0$ because it is not a saddle. 
Instead, we can use the unique 2-cycle $\{Q, F(Q)\}$
which remains of saddle type throughout the bifurcation.

The left hand side of Figure 2a shows this for $a=.3$; the point $Q$ corresponds to
the tip at the rightmost point, and the factor for self-similarity is approximately 8.
The two pictures, Figure 1 left and Figure 2a left, are slices at different
points $O$ and $Q$ of  $\partial \cB$. However, the ``tip'' shape of the slice
$W^u(Q)\cap \cB$ at $Q$ appears to be repeated densely at small scales in the
slice $T\cap \cB$ as well as in $W^u(Q)\cap \cB$.
This might be explained by the existence of transversal intersections between the
stable manifold $W^s(Q)$ and $T$ at a dense subset of $T\cap \partial \cB$
Similarly the ``cusp'' at $O$ of the slice $T\cap \cB$ appears to be repeated
at small scales in the slice $W^u(Q)\cap \cB$ 
as well as in $T\cap \cB$.
A similar result, showing that all unstable slices have features in common has been 
formulated and proved for hyperbolic maps in [BS7].
 
In general, the set $K^+$ where the orbits are bounded, coincides with $\{G^+=0\}$.  This set is of primary interest, even if $\cB=\emptyset$ .  In Figures 1 and 2, $\cB$ seems to have ``exploded''  leaving $K^+ = \partial K^+ = J^+$ without interior   when $\epsilon \neq 0$.  On the other hand, the computer detail in Figure 2b persists as $\epsilon \to 0$. This means that $\epsilon \mapsto J^+(F_\epsilon)$ will appear to have  ``exploded"  a little bit as $\epsilon \neq 0$.  We will see a marked similarity between Figures~2b and 5, which corresponds to Theorems 1 and 2.

We may use $\varphi^{\iota,o}$ to represent ${\cal B}\cap\Sigma$ graphically for the H\'enon family $F(x,y)$ discussed above.  We may use  $\varphi^o$ to parametrize $\Sigma$.  In Figure 3 we have drawn  part of the slice $\cB\cap \Sigma$,   with level sets of the real and imaginary parts of $\varphi^\iota$.  One critical point for $h_\alpha$ (as well as its complex conjugate and translates) is clearly evident on the left hand picture, and at least two more critical points are evident on the right.   Figures 1 and 3 give invariant slices of the same basins and share certain features, but Figure 3, which is specialized to parabolic basins, has more focus on the interior; and the two pictures are localized differently.

\section{``Almost Fatou'' coordinates: dimension 1.}   Consider a family of maps
$$f_{\epsilon} (x) = x + (x^2+\epsilon^2)\alpha_\epsilon(x), \ \ \alpha_\epsilon(x) = 1 + p\,\epsilon + (q+1)x + \cdots. \eqno(2.1)$$
  We are interested in analyzing $f^n_\epsilon (x)$ for $n$ large and $\epsilon$ small and $f^n_\epsilon (x)$ near 0.  The first step is to introduce a change of coordinates in which $f_\epsilon$ is close to a translation.    If we change coordinates to $(\hat x,\hat\epsilon)$ given by  $x = (1-p\,\hat\epsilon)\hat x$, with $\epsilon= \hat\epsilon - p\,\hat\epsilon^2$, then we have $p=0$.
Define
$$\gamma_\epsilon(x) = {\alpha_\epsilon(x)\over 1+ x\alpha_\epsilon(x)} = 1 +  q x + \cdots  \eqno(2.2)$$
Let  $\epsilon\in{\bf C}$ be such that 
$$ 0<\Re\epsilon, \ \ \ |\Im\epsilon|\le  {\rm const. }\ |\epsilon|^2   \eqno(2.3)$$
We consider the coordinate change 
$$x\mapsto {u}_\epsilon ={1\over \epsilon}\arctan {x\over\epsilon} =  {1\over 2i\epsilon} \log {i\epsilon -x\over i\epsilon + x} \eqno(2.4)$$

\bigskip

\centerline{\includegraphics[width=4.1in]{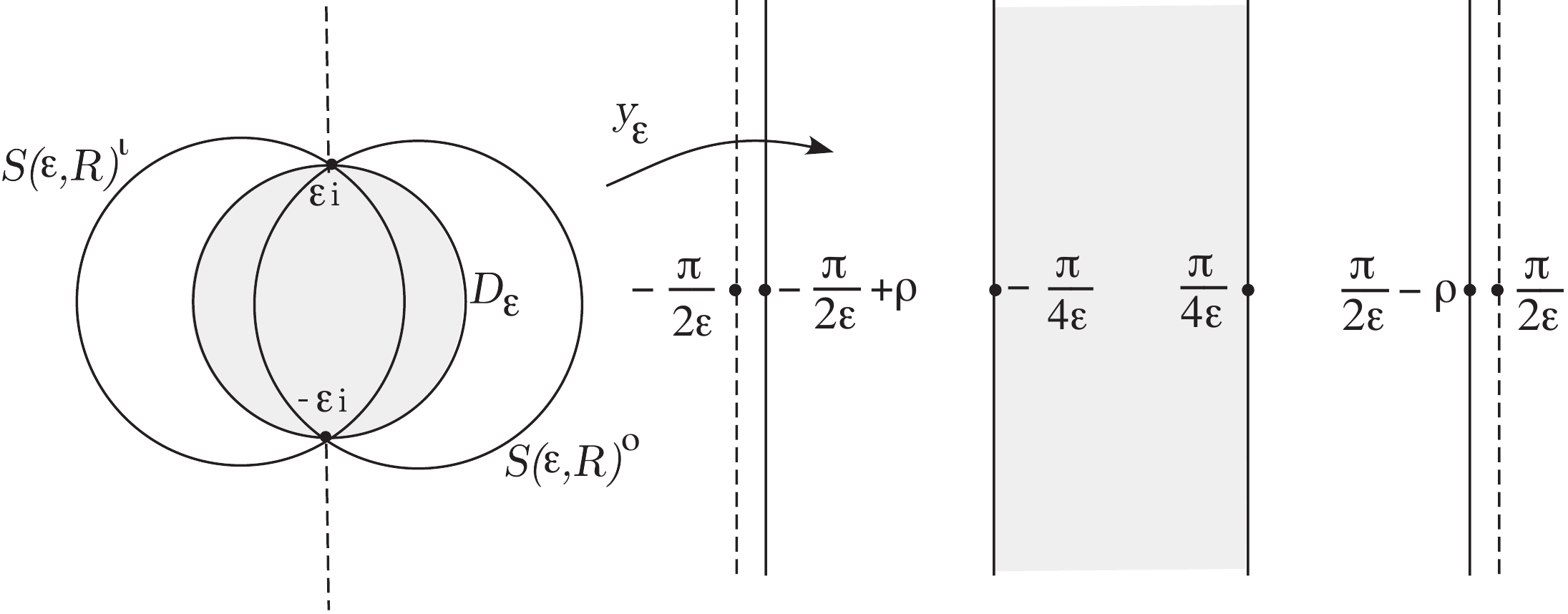} }

\centerline{Figure 4.  Mapping of the slit region by ${u}_\epsilon$ for $\epsilon>0$.    } 
\bigskip

To describe this coordinate change, we let $S(\epsilon,R)$ be the union of two disks of radius $R$ as in Figure 4.   $R$  will be chosen small enough that the $\cdots$ terms in (2.1) and (2.2) are small.  Note, however, that the proportions in Figure 4 may be misleading because $R$ will be fixed while $\epsilon\to0$, so $S(\epsilon,R)$ will be of a fixed diameter while $D(\epsilon)$ shrinks.  Let $H_\epsilon^\iota $ denote the half-space to the left of the line $\epsilon i{\bf R}$; this is the space of ``incoming'' points.  Since $\epsilon$ is almost real, $H^\iota_\epsilon$ is approximately the left half plane.   Figure 4 shows the boundary of $H^\iota_\epsilon$ (dashed) and its image (also dashed).   We let $S(\epsilon,R)^\iota=S(\epsilon,R)\cap H_\epsilon^\iota$ denote the incoming half of $S(\epsilon,R)$.  The image ${u}_\epsilon(S(\epsilon,R)^\iota)$ is bounded by of two parallel lines: one of them passes through $-{\pi\over 2\epsilon}+\rho$, and the other is inside the shaded strip.   Similarly, we define $H_\epsilon^o$ to be the  half-plane to the right, and $S(\epsilon,R)^o=S(\epsilon,R)\cap H_\epsilon^o$ to be the ``outgoing'' points.  The image of $D_\epsilon:=\{|x|<|\epsilon|\}$ (shaded) is the shaded strip on the right hand side.  An important feature of this picture is that if $R$ is fixed, then $\rho$ stays bounded as $\epsilon\to0$.

From (2.1) and (2.2) we have
$$ {i\epsilon-f_\epsilon(x)\over i \epsilon + f_\epsilon(x)} = {(i\epsilon-x) \{ 1+(x+i\epsilon)\alpha_\epsilon(x) \} \over (i\epsilon + x) \{ 1+(x-i\epsilon)\alpha_\epsilon(x) \}} =  { (i\epsilon-x) (1+ i\epsilon\gamma_\epsilon (x)) \over (i\epsilon+x)(1- i\epsilon\gamma_\epsilon(x))}$$
Thus we have
\begin{equation}
\begin{split}
{u}_\epsilon(f_\epsilon (x))-{u}_\epsilon (x)  &   = {1\over 2i\epsilon}  \log {1+i\epsilon\gamma_\epsilon(x) \over 1 - i\epsilon\gamma_\epsilon(x)} \cr
 & = {1\over i\epsilon} \left\{i\epsilon\gamma_\epsilon(x) + {1\over 3} (i\epsilon\gamma_\epsilon(x))^3 +\cdots\right\}\cr
& =  \gamma_\epsilon(x) -{\epsilon^2\over 3}\gamma_\epsilon(x)^3 +\cdots \cr
& =    1+ qx + O(|\epsilon|^2 + |x|^2) 
\end{split}
\tag{2.5}
\end{equation}

\begin{proposition}
For any compact subset $C\subset S^\iota(R)$, there are $\epsilon_0>0$ and $C_0, K_0>0$ such that for $|\epsilon|<\epsilon_0$ and $x\in C$, the following hold:
\item{(i)}  $f^j_\epsilon(x)\in S^\iota(\epsilon,R)\cup D_\epsilon$, for $0\le j\le {3\pi\over 5|\epsilon|}-K_0$
\item{(ii)} $|f^j_\epsilon(x)|\le C_0 \max \left\{ {2\over j},|\epsilon|\right \}$, for $0\le j\le {3\pi\over 5|\epsilon|}-K_0$
\item{(iii)} $f^j(x)\in D_\epsilon$ for ${\pi\over 3|\epsilon|}\le j\le {3\pi\over 5|\epsilon|}-K_0$. 
\end{proposition}

\noindent{\it Proof.}   First we note that
$${u}_\epsilon(x) + {\pi\over 2\epsilon} \to -{1\over x} \ \ \ (\epsilon\to0)$$
uniformly on compact subsets of $S^\iota(0,R)$.  So there is a $K_0>0$ such that
$$-{\pi\over 2|\epsilon|} < \Re \left({\epsilon\over|\epsilon|} {u}_\epsilon(x)\right) < -{\pi\over 2|\epsilon|}+K_0$$
on $C$.  We have
$${3\over 4} < \Re\left({\epsilon\over|\epsilon|} {u}_\epsilon(f_\epsilon(x))\right) - \Re\left({\epsilon\over|\epsilon|}{u}_\epsilon(x)\right) < {5\over 4}$$
and so it follows that
$$-{\pi\over 2|\epsilon|} + {3j\over 4} < \Re\left( {\epsilon\over|\epsilon|}{u}_\epsilon(f_\epsilon^j(x))\right) < -{\pi\over 2|\epsilon|} + {5j\over 4} + K_0.  \eqno(2.6)$$
If $0\le j<3\pi/5|\epsilon|-K_0$, then
$$\Re\left({\epsilon\over|\epsilon|}{u}_\epsilon(f^j_\epsilon(x))\right)<{\pi\over 4|\epsilon|}$$
and hence $f^j_\epsilon(x)\in S^\iota(\epsilon,R)\cup D_\epsilon$, which proves $(i)$.

For $(ii)$ we note that the function $\tan z$ maps any line $\Re z = a$ ($a\in {\bf R}$) to the circular arc with endpoints $\pm i$ and passing through the point $\tan a$.  It follows that 
$|\Re z|<{\pi\over 2}\Rightarrow |\tan z|\le 1$, and ${\pi\over 4}\le |\Re z|<{\pi\over 2}$ implies that $|\tan z| \le \tan|\Re z|<{1\over {\pi\over 2}-|\Re z|}$.  By (2.4) we have
$$-{\pi\over 2} \le \Re(\epsilon {u}_\epsilon) \le -{\pi\over 4} \Rightarrow |x|\le |\epsilon| \tan|\Re(\epsilon {u}_\epsilon)|<{ |\epsilon|\over {\pi\over 2} + \Re(\epsilon {u}_\epsilon)}$$
and $|\Re (\epsilon {u}_\epsilon)|\le {\pi\over 4}$ implies that $|x|\le |\epsilon|$.  Now by (2.6) we have
$${\pi\over 2|\epsilon|} + \Re\left( {\epsilon\over |\epsilon|}{u}_\epsilon (f_\epsilon(x)\right ) \ge cj$$
for some $c$.  This shows $(ii)$.

If $\pi/(3|\epsilon|)\le j < 3\pi/(5|\epsilon|) - K_0$, then by (2.6)
$$-{\pi\over 4|\epsilon|}<\Re\left({\epsilon\over|\epsilon|}{u}_\epsilon(f^j_\epsilon(x)) \right ) < {\pi\over 4|\epsilon|}$$
and hence $f^j_\epsilon(x)\in D_\epsilon$, which proves $(iii)$.  \qed

Next, with $q$ as in (2.2), we define
\begin{align*}
w_\epsilon(x) &  =  {u}_\epsilon(x) - {q\over 2} \log \left( 1 + {x^2\over \epsilon^2} \right) \cr
& = {1\over 2i\epsilon}  \log { i\epsilon -x\over i\epsilon + x} - {q\over 2} \log(\epsilon^2 + x^2) + q\log\epsilon 
\end{align*}
The corresponding incoming and outgoing versions are obtained by adding terms that depend on $\epsilon$ but do not depend on  $x$:
\begin{align*}
w_\epsilon^{\iota/o} : & = w_\epsilon(x) - q \log\epsilon  \pm {\pi\over 2\epsilon} \cr
& = {1\over\epsilon} \left(  \pm {\pi\over 2}+\arctan {x\over \epsilon} \right ) -{q\over 2}  \log(\epsilon^2+x^2)
\end{align*}

\begin{lemma}
$\lim_{\epsilon\to0} w_\epsilon^\iota = w_0^\iota$, and $\lim_{\epsilon\to0} w_\epsilon^o = w_0^o$.
\end{lemma}

Let us define $A_\epsilon(x):=w_\epsilon(f_\epsilon(x))-w_\epsilon(x)-1$, which measures how far $w_\epsilon(x)$ is from being a Fatou coordinate. 

\begin{proposition}
$A_\epsilon(x) = O(|\epsilon|^2+|x|^2)$.
\end{proposition}

\noindent{\it Proof.}  First we observe that
\begin{align*}
\epsilon^2 + f_\epsilon(x)^2 &  = \epsilon^2 + \{x + (\epsilon^2 + x^2) \alpha_\epsilon(x)\}^2\cr
& = \epsilon^2+x^2 + 2x(\epsilon^2 + x^2)\alpha_\epsilon(x) + (\epsilon^2+x^2) \alpha_\epsilon(x)^2 \cr
& = (\epsilon^2 + x^2) \left( 1 + 2x \alpha_\epsilon(x) + (\epsilon^2+x^2)\alpha_\epsilon(x)^2 \right)
\end{align*}
Thus
\begin{align*}
\log \left(1 + {f_\epsilon(x)^2\over \epsilon^2} \right) - \log \left( 1 + {x^2\over\epsilon^2}\right) & = \log(1 + 2x\alpha_\epsilon(x)+ \cdots)\cr
& = 2x\alpha_\epsilon(x) + O(x^2) \cr
& = 2x + O(|(\epsilon,x)|^2)
\end{align*}
It follows that 
\begin{align*}
w_\epsilon(f_\epsilon & (x))  -w_\epsilon(x) = \cr
& = ({u}_\epsilon(f_\epsilon(x))-{u}_\epsilon(x)) + {q\over 2}\left(  \log\left( 1 + {f_\epsilon^2(x)\over \epsilon^2}\right) - \log\left( 1 + {x^2\over \epsilon^2}\right) \right)\cr
& = (1 + qx + O(|\epsilon|^2+|x|^2)  - {q\over 2}\left( 2x + O(|\epsilon|^2+|x|^2) \right)\cr
& = 1+ O(|\epsilon|^2+|x|^2) \cr 
\end{align*}
which gives the desired result.  \qed
We note, too, that
$$A_\epsilon(x) = A_0(x) + \epsilon \tilde A(x) + O(\epsilon^2)  \eqno(2.7)$$
where $\tilde A(x) = O(x)$. 

\begin{corollary}
$w^{\iota/o}_\epsilon(f_\epsilon(x))-w^{\iota/o}_\epsilon(x) - 1 =  O(|\epsilon|^2+|x|^2)$
\end{corollary}

\begin{lemma}
There exists $K_0>0$ such that:  If $x,f_\epsilon(x),\cdots,f_\epsilon^n(x)\in S(\epsilon,R)$, then
$$|w_\epsilon(f^n_\epsilon(x))-w_\epsilon(x)-n|\le K_0$$
and hence
$$|w^o_\epsilon(f^n_\epsilon(x)) - w^\iota_\epsilon(x)+ {\pi\over\epsilon}-n|\le K_0$$
\end{lemma}

\noindent{\it Proof.}  We have 
$$w_\epsilon(f^n_\epsilon(x))- w_\epsilon(x)-n = \sum_{j=0}^{n-1} A_\epsilon(f^j_\epsilon(x)).$$
Choose $0<n_1<n_2<n$ such that
\begin{align*}
 & f^j_\epsilon(x) \in S^\iota(\epsilon,R)-D_\epsilon, \ \ \  0\le j\le n_1-1 \cr
 & f^j_\epsilon(x) \in D_\epsilon, \ \ \ \ \ \ \ \ \ \ \ \ \ \ \ \ \  n_1\le j\le n_2-1 \cr
 &f^j_\epsilon(x)\in  S^o(\epsilon,R)-D_\epsilon, \ \ \ n_2\le j\le n
\end{align*}
Then $n_2-n_1\le const/|\epsilon|$, and
\begin{align*}
&  |A_\epsilon(f^j_\epsilon(x))| \le const/j^2, \ \ \ 0\le j\le n_1-1 \cr
&  |A_\epsilon(f^j_\epsilon(x))| \le const \, |\epsilon|^2, \ \ \ n_1\le j\le n_2-1\cr
&  |A_\epsilon(f^j_\epsilon(x))| \le const/(n-j)^2, \ \ \ n_2\le j\le n
\end{align*}
This proves the Lemma.  \qed

We will use the following condition:
$$\{n_i,\epsilon_i\} {\rm\ is\ a\ sequence\ such\ that\ \  }{\pi\over 2\epsilon_i} -n_i{\rm\ \   is\ bounded}  \eqno(2.8)$$
Recall that $\{j_i,\epsilon_i\}$ is an {\it  $\alpha$-sequence} if $\epsilon_i\to0$, and  $j_i-{\pi\over  \epsilon_i}\to \alpha$ as $i\to\infty$.  For instance, $(m,\epsilon_m)$ with $\epsilon_m ={\pi\over m-\alpha}$ is an $\alpha$ sequence.   Every $\alpha$-sequence satisfies (2.3), and if $\{j_i,\epsilon_i\}$ is an $\alpha$-sequence, then $\{j_i/2,\epsilon_i\}$ satisfies (2.8)

\begin{proposition}
If (2.8) holds, and if $C$ is a compact subset of $S^\iota(R)$, then $\{f^{n_i}_{\epsilon_i}\}$ is uniformly bounded and forms a normal family on $C$.
\end{proposition}

We define an almost Fatou coordinate in the incoming direction:
$$\varphi^\iota_{\epsilon, n}(x) = w^\iota_\epsilon(f^n_\epsilon(x)) -n = w^\iota_\epsilon(x) +\sum_{j=0}^{n-1}A_\epsilon(f^j_\epsilon(x)).$$
We recall that $\cB$ denotes the parabolic basin of points where the iterates $f^j$ converge locally uniformly to  $O = (0,0)$.

\begin{theorem}
If  (2.8) holds,  then on $\cB$ we have
$$\lim_{j\to\infty}\varphi^\iota_{\epsilon_j,n_j}=\varphi^\iota.$$   
\end{theorem}

\noindent{\it Proof.}  If $x\in\cB$, we may assume that $x\in S(\epsilon,R)^\iota$, where $\epsilon$ and $R$ are as above.  If we set $\varphi^\iota_{0,n}=w_0^\iota + \sum_{j=0}^{n-1} A_0(f^j_0 x)$, we have $\varphi^\iota=\lim_{n\to\infty}\varphi^\iota_{0,n}$.  We consider
 $$\varphi^\iota_{\epsilon,n}-\varphi^\iota_{0,n}= w^\iota_\epsilon(x) -w^\iota_0(x) + \sum_{j=0}^{n-1}\left (A_\epsilon(f^j_\epsilon(x)) - A_0(f^j_0(x))\right).$$
 We will show that this difference vanishes as $\epsilon=\epsilon_j\to0$ and $n=n_j\to\infty$.  We have $w^\iota_\epsilon-w^\iota_0\to0$ by Lemma 2.2.  The summation is estimated by
 $$\left| \sum\right| \ \le \sum\left| A_0(f^j_\epsilon(x))-A_0(f_0^j(x)) \right| + \sum \left| A_\epsilon(f^j_\epsilon(x)) - A_0(f^j_\epsilon(x))\right| = {\sum}_1 + {\sum}_2$$
For the first sum, we recall that $A_0(x)=O(x^2)$, and so by Proposition 2.1 we have that the two series are summable:
$$\sum_{j=1}^{n-1} \left|A_0(f^j_\epsilon)\right| +  \left|A_0(f^j_0)\right| \le  K \sum_{j=1}^{{\pi\over 2|\epsilon|}}\left( {1\over j^2} + |\epsilon|^2 \right)\le K\left( {\pi|\epsilon|\over 2} +\sum_{j=1}^\infty{1\over j^2} \right)\le  B$$
as $n\to\infty$ and $\epsilon\to0$.  For $\delta>0$ we choose  $J$ such that $\sum_J^\infty j^{-2}<\delta$.   If we write ${\sum}_1=\sum_1^J + \sum_{J+1}^\infty$, then we see that $\sum_{J+1}^\infty\le \pi K|\epsilon|/2 + \delta$.  On the other hand,  for fixed $j$ we have $A_0(f^j_\epsilon)\to A_0(f^j_0)$ as $\epsilon\to0$, so we conclude that 
$$\sum_{1}^{J} = \sum_{1}^{J} \left |A_0(f^j_\epsilon)-A_0(f^j_0) \right|\to0$$ 
as $\epsilon\to0$.  In conclusion, we see that $\lim_{\epsilon\to0}{\sum}_1\le K\delta$ for all $\delta$, so that ${\sum}_1\to0$.

For the second part, we use (2.7) so that
\begin{align*}
{\sum}_2 &\le  \sum_{j=0}^{n-1} \left| A_\epsilon(f_\epsilon^j)-A_0(f_\epsilon^j)\right| =  \sum_{j=0}^{n-1} \left| \epsilon \tilde A(f^j_\epsilon)\right|+ \sum_{j=0}^{{\pi\over 2|\epsilon|}}K|\epsilon|^2 \cr
& \le K'|\epsilon| +K' |\epsilon|  \sum_{j=1}^{{\pi\over 2|\epsilon|}} {1\over j}  \le  K''|\epsilon|  \log\left( {\pi\over 2|\epsilon|}\right)
\end{align*}
and this last term vanishes as $\epsilon\to0$, which completes the proof.
 \qed

Using $f^{-1}_\epsilon$ and $w^o_\epsilon$, we may also define almost Fatou coordinates in the outgoing direction, and the direct analogue of Theorem 2.7 holds:

\begin{corollary}
The inverse maps $(\varphi^o_{\epsilon_j,n_j})^{-1}$ converge uniformly to $(\varphi^o)^{-1}$ on compact subsets of  $\varphi^o(S(R)^o)$.
\end{corollary}

Let us consider $T_\alpha:= (\varphi^o)^{-1}\circ\tau_\alpha\circ\varphi^\iota$.  For  $\alpha\in\bC$, we define $D_\alpha\subset\bC$ to be the set where $T_\alpha$ is defined.  
The quantities $\varphi^{\iota/o}-w^{\iota/o}$ are bounded, so the range of $\varphi^\iota$ is approximately $\{\Re(\zeta)<K\}$ and the range of $\varphi^o$ is approximately $\{\Re(\zeta)<-K\}$.   Thus, for each $x\in S(R)^\iota$, we have $D_\alpha\ne\emptyset$ if $\Re(\alpha)$ is sufficiently negative.

From Theorem 2.7 we have:

\begin{theorem}
If $(\epsilon_j,n_j)$ is an $\alpha$-sequence, then $\lim_{j\to\infty}f_{\epsilon_j}^{n_j}=T_{\alpha}$  on $D_\alpha$. 
\end{theorem}

\noindent{\it Proof.}    We may assume that $x\in S(\epsilon,R)^\iota$.  Write $x'=f^{n_j}_{\epsilon_j}(x)$.  Choose $m_j$ and $m_j'$ so that $n_j=m_j+m'_j$, and (2.8) holds for $\{m_j\}$ and $\{m_j'\}$.   Thus  $f^{m_j}(x) = f^{-m_j'}(x'_j)$, so $w_{\epsilon_j} f^{m_j}(x) = w_{\epsilon_j} f^{-m_j'}(x'_j)$.  We rewrite this as
$$\varphi^\iota_{\epsilon_j,m_j}(x) + m_j -{\pi\over 2 \epsilon_j} = \varphi^o_{\epsilon_j,m_j'}(x'_j) - m_j' +{\pi\over 2\epsilon_j},$$
or
$$\varphi^\iota_{\epsilon_j,m_j}(x) + n_j-{\pi\over\epsilon_j}  = \varphi^o_{\epsilon_j,m'_j}(x'_j).\eqno(2.9)$$

By Theorem 2.7, $\varphi^\iota_{\epsilon_j,m_j}(x)$ converges to  $\varphi^\iota(x)= \zeta_0\in\bC$.  Since $x\in D_\alpha$, we know that $\zeta_0+\alpha$ is in the range $\varphi^o(S(R)^o)$.   If we replace $x$ by a preimage $f^{-k}x$, we will have $\varphi^\iota(x)=\zeta_0-k$.   We may assume that $\zeta_0-k$ is in the range of $\varphi^o$, and so $(\varphi_{\epsilon_j,m'_j}^o)^{-1}$ converges to $(\varphi^o)^{-1}$ in a neighborhood of $\zeta_0-k$.  It follows that the points $x_j'=f_{\epsilon_j}^k(\varphi^o_{\epsilon_j,m'_j})^{-1}(\zeta_0-k)$ converge to a limit $x'$.  As $j\to\infty$, we may pass to a limit in (2.9) to obtain $\varphi^\iota(x) + \alpha =\varphi^o(x')$.   Applying $(\varphi^o)^{-1}$ to both sides of the equation, we see that $T_{\alpha}x=x'$.  \qed

\medskip

\section{Two-dimensional case: Convergence of the ``Almost Fatou'' coordinate.}
We consider a one-parameter family $F_\epsilon$, varying analytically in $\epsilon$, such that  $$F_0(x,y) = (x+x^2+\cdots, by + \cdots).$$  The fixed point $O=(0,0)$ has multiplicity 2 as a solution of the fixed point equation, and we will assume that for $\epsilon\ne0$ the fixed point $O$ will split into a pair of fixed points.    We parametrize so that the fixed points are $(\pm i\epsilon,0)+O(\epsilon^2)$.   We consider here only fixed points of multiplicity two.  We suspect that perturbations of fixed points of higher multiplicity might be quite complicated, since this is  already the case in dimension 1, as was shown by Oudkerk [O1,2].

\begin{theorem}
By changing coordinates and reparametrizing $\epsilon$, we may suppose  that our family of maps has the local form
$$F_\epsilon(x,y) =   \left(x + (x^2+\epsilon^2)\alpha_\epsilon(x,y), b_\epsilon(x) y + (x^2+\epsilon^2)\beta_\epsilon(x,y)  \right)   \eqno(3.1)  $$
where $\alpha_\epsilon = 1 + (q+1)x + ry + O(|x|^2+|y|^2+|\epsilon|^2)$, $b_0(0)=b$.  In particular, the points $(\pm i\epsilon, 0)$ are fixed, the lines $\{x=\pm i\epsilon\}$ are local stable manifolds, and the map is locally linear on the stable manifolds.  Further, the multipliers at the fixed points are  $(1\pm 2i\epsilon, b_\epsilon(\pm i\epsilon))+O(\epsilon^2)$. 
\end{theorem}

\noindent{\it Proof.}  By a change of variables, we may assume that the fixed points are $(\pm i \epsilon,0)$.  Each fixed point will have eigenvalues $1+O(\epsilon)$ and $b+O(\epsilon)$.  There will be local strong stable manifolds corresponding to the eigenvalue $b+O(\epsilon)$.  We apply the graph transform  as in [HPS, \S5A]  in order to obtain a domain for the stable manifold which is uniformly large in $\epsilon$.  Rescaling coordinates, we may assume that we have graphs $$W^{s}_{loc}(\pm i\epsilon,0)=\{x=\psi^\pm(\epsilon,y):|y|<1\},$$
where $\psi^\pm$ is analytic in $\epsilon$ and $y$.  Let us consider new coordinates $X,Y$ defined by  $x= \chi_0(\epsilon,y)+X\chi_1(\epsilon,y)$, $Y=y$, where we set $\chi_0=-{1\over 2}(\psi^++\psi^-)$ and $\chi_1={1\over 2 i \epsilon}(-\psi^++\psi^-)$.  Since  $\psi^\pm$ are uniquely determined and analytic in $\epsilon$, we have $\lim_{\epsilon\to0}(\psi^+(\epsilon,y)-\psi^-(\epsilon,y))=0$, from which we conclude that $\chi_1$ is analytic in $(\epsilon,y)$. 

In order for $F$ to have the desired form in the $y$-coordinate, we need to change coordinates  so that   $y\mapsto F_\epsilon(\pm i\epsilon, y)$ is linear in $y$.  We set $b_\epsilon^\pm={\partial F\over \partial y}(\pm i\epsilon, 0)$.  There is a unique function $\xi^\pm_\epsilon(y)= y + O(y^2)$ such that $F(\pm i\epsilon,\xi_\epsilon^\pm(y)) = b^\pm_\epsilon \xi_\epsilon^\pm(y)$.  We note that $\xi^\pm_\epsilon$ is holomorphic in $\epsilon$, and $\xi^+_0=b\xi^-_0$.  Thus $\epsilon\mapsto (\xi^-_\epsilon -\xi^+_\epsilon)/\epsilon$ is analytic, and we may define a new coordinate system $(X,Y)$ with $X=x$ and $Y = [(i\epsilon-x)\xi^-_\epsilon(y) + (x+i\epsilon)\xi^+_\epsilon(y)]/(2i\epsilon)$.  $F$ has the desired form in the new coordinate system.  

Our map now has the form (3.1) with $\alpha_\epsilon = 1 + p\,\epsilon + (q+1)x + ry +O(|x|^2+|y|^2+|\epsilon|^2)$.  We can make $p=0$ using the coordinate change   $x = (1-p\,\hat\epsilon)\hat x$, with $\hat\epsilon$ defined by $\epsilon= \hat\epsilon - p\,\hat\epsilon^2$.  The remaining statements in the Theorem are easy consequences of (3.1).  \qed

One motivation for the normalization in (3.1) is that for the map $z\mapsto z+z^2+\epsilon^2$, the fixed points are $\pm i\epsilon$, and the multipliers are $1+2i\epsilon$.

\bigskip\noindent{\bf Remarks about $\alpha$-sequences.}   If we wish to use Theorems 1 and 2 for a specific family of mappings, we need first to make the changes of coordinates involved with Theorem 3.1.  This influences the value of $\alpha$ which appears in the $\alpha$-sequence $\{\epsilon_j\}$.  That is, if $\epsilon_j$ is an $\alpha$-sequence, and if $\epsilon=\hat\epsilon-p\hat\epsilon^2$, then $\hat\epsilon_j$ is an $\hat\alpha$-sequence with $\hat\alpha = \alpha+\pi p$.  

Let us discuss how the normalization relates to the H\'enon maps $F_{a,\epsilon}$ in (1.1).  If we first make a linear change of coordinates so that the axes are the eigen-directions of $Df_\epsilon$ at $O$, then $F_{a,\epsilon}$ becomes  $(x+ x^2/(1-a) +\epsilon^2 + a y (2x+ ay)/(1-a), a y + O(|x|^2+|y|^2) )$; and in particular the $p\,\epsilon$ term, discussed above, vanishes.    Now in order to have the form (3.1) we conjugate with the dilation $x\mapsto (1-a)x$.  This gives us $(x+ x^2 +\epsilon^2 / (1-a),  a y)+\cdots$.  Thus for condition (3.1) we use the parameter $\epsilon'=\epsilon/(1-a)^{1/2}$, and this means that for Theorems 1, 2, and 3.9, it is $\epsilon'$ which must be part of an $\alpha$-sequence.

\bigskip
As in the 1-dimensional case we define
$$\gamma_\epsilon(x,y) = {\alpha_\epsilon(x,y)\over 1 + \alpha_\epsilon(x,y)} = 1 + qx + ry + \cdots.$$
If we define $\tilde y_\epsilon(x,y) = y_\epsilon(x)$, then as in \S2 we will have
\begin{equation}
\begin{split}
\tilde y_\epsilon(F_\epsilon(x,y))-\tilde y_\epsilon(x,y) & =  \gamma_\epsilon(x) -{\epsilon^2\over 3}\gamma_\epsilon(x)^3 +\cdots \cr
& =    1+ qx + ry+ O(|\epsilon|^2 + |x|^2+|y|^2)
\end{split}\tag{3.2}
\end{equation}
We use the notation $(x_j,y_j):= F^j_\epsilon(x,y)$, and we let $R$ and $\epsilon$ be as in \S2.  For $\eta_0$, we define $\tilde  S(\epsilon,R)^\iota=  S(\epsilon,R)^\iota\times\{|y|<\eta_0\}$.   Arguing as in Proposition 2.3, we have

\begin{proposition}
We may choose $\eta_0>0$ small enough that if  $|\epsilon|<\epsilon_0$, and $C\subset \tilde S(\epsilon,R)^\iota$ is compact, then there exists $K<\infty$ such that $|x_{\epsilon,j}|\le K\max({1\over j}, |\epsilon|,|b|^n)$ for $j\le {\pi\over 2|\epsilon|}$ and $(x,y)\in C$.
\end{proposition}

Following \S2, we define 
$$\tilde w_\epsilon(x,y)=w_\epsilon(x) + ry/(b-1), \ \ \ \tilde w^{\iota/o}_\epsilon(x,y)=w^{\iota,o}_\epsilon(x) + ry/(b-1).$$  
As in Proposition 2.3, we have:

\begin{proposition}
$\tilde A_\epsilon:= \tilde w_\epsilon(F_\epsilon(x,y)) -\tilde w_\epsilon(x,y) = 1 + O(|\epsilon|^2 + |x|^2 + |y|^2)$.
\end{proposition}

We define the incoming almost Fatou coordinate:
$$\varphi^\iota_{\epsilon,n}(x,y): = \tilde w^\iota_\epsilon(F^{n}_\epsilon(x,y)) - n$$
and as in \S2, we obtain:

\begin{theorem}
If $\epsilon_j\to 0$, and if $n_j$ satisfies (2.8), then $\lim_{j\to\infty}\varphi^\iota_{\epsilon_j,n_j} = \varphi^\iota_0$ locally uniformly on $\cB$.
\end{theorem}

We omit the proof since it is essentially the same as in the 1-dimensional case.  By the Center Stable Manifold Theorem, there is a 1-parameter family of  manifolds $W^{cu}_{\epsilon, loc}$, $|\epsilon|<\epsilon_0$, the (local) center unstable manifolds of the fixed points of $F_\epsilon$, corresponding to the (larger) eigenvalue, which is near 1.  We denote them by $M_\epsilon$ and note that there is a neighborhood $U$ of $O$, with the property that $U\cap f_\epsilon(M_\epsilon)\subset M_\epsilon$.  The manifolds $M_\epsilon$ can be taken to be $C^1$ smooth and to vary in a $C^1$ fashion with respect to $\epsilon$.

\begin{proposition}
Let us fix a compact $W\subset U$.  There are constants $\beta<1$ and $C<\infty$ such that for each $p\in W$,  ${\rm dist}(F^j_\epsilon p,M_\epsilon)\le C\beta^j$ for $1\le j\le j_0$ if $f^j_\epsilon p \in U$ for $1\le j\le j_0$.  
\end{proposition}

Let $\pi(x,y)=x$ be the projection to the $x$-axis.  The tangent space to $M_\epsilon$ is close to the $x$-axis, so there is a function $h_\epsilon$ such that $M_\epsilon=\{(x,h_\epsilon(x)):|x|<R\}$.    
By [U2] there is an analytic function $h$ on $\{|x-r|<r\}$ which extends continuously to the closure, and which satisfies $h(0)=0$, and $\Sigma_0:=\{(x,y):|x-r|<r, y=h(x)\}$ is contained in $\Sigma$.  Further, $\Sigma_0$ is invariant in negative time, so $h$ coincides with the function $h_0$, and $\Sigma_0\subset M_0$. 

\begin{proposition}
With the hypotheses of Proposition 3.5, let us suppose that $p\in W$ and $F^jp\in U$ for $1\le j\le j_0$.  Then  ${\rm dist}(F^jp,\Sigma_0) \le C(|\epsilon|+\beta^j)$ for  those values of $j$ for which $1\le j\le j_0$, and $\pi F^j_\epsilon p\in\{|x -r|<r \}$.
\end{proposition}

Let $T(\epsilon,-n)=\{p\in U: F_\epsilon ^{-j}p\in U: 0\le j\le n\}$.  We define an outgoing almost Fatou coordinate for $p\in T(\epsilon,-n)\cap M_\epsilon$ by setting
$$\varphi^o_{\epsilon,n}(p) =w^o_\epsilon(F^{-n}p) + n$$
\begin{proposition}
Suppose that the sequence $(\epsilon_j,n_j)$ satisfies (2.8) and that $|x-r|<r$.  Then
$$\lim_{j\to\infty} \varphi^o_{\epsilon_j,n_j}(x,h_{\epsilon_j}x) = \varphi^0(x,h_0x).$$
\end{proposition}

\noindent{\it Proof.}  By \S2, we have that $p$ will belong to $T(\epsilon_j,-n_j)$ for $j$ sufficiently large.  The proof of this Proposition is then essentially the same as the proof of Theorem 2.7.
\qed

\begin{proposition}
Suppose that $p\in M_\epsilon$ and that the projection $x=\pi_1( p)$ satisfies $|x+r|<r$.  Suppose, further, that $T_\alpha p\in\Sigma_0$.   If $\{\epsilon_j\}$ is an $\alpha$-sequence, then $F^{n_j}_{\epsilon_j}p$ converges to $T_\alpha p$.
\end{proposition}

\noindent{\it Proof.}  Shrinking $U$ if necessary, we may choose $\beta$ and $\hat\beta$ such that $\beta<1$, $\beta^2\hat\beta<1$ with the following properties:  in the vertical direction, $F$ contracts with a factor of $\beta$;  and ${\rm dist}(F^{-1}q_1,F^{-1}q_2)\le \hat\beta\, {\rm dist}(q_1,q_2)$.   Now let us write $q=F^{n_j}_{\epsilon_j}p$, and let $q':=(\pi q,h_{\epsilon_j}(\pi q))$ denote the projection to $M_{\epsilon_j}$.  By Proposition 3.6, we have ${\rm dist}(q,q')=O(\beta^{n_j})$.  Now we write $n_j=m_j'+m_j''$, where $m_j'$ and $m_j''$ are both essentially $n_j/2$.  We have 
$$F^{m_j'}_{\epsilon_j}p = F^{-m_j''}_{\epsilon_j}q = F^{-m''_j}_{\epsilon_j}q' + \hat\beta^{m_j''}\beta^{n_j} =  F^{-m''_j}_{\epsilon_j}q' + o(1)$$
Adding and subtracting $\pi/(2\epsilon_j)$ and $n_j$ to  $w_{\epsilon_j}F^{m_j'}_{\epsilon_j}p = w_{\epsilon_j} F^{-m_j''}q' + o(1)$, we have
$$w^{\iota}_{\epsilon_j}F^{m_j'}p - m_j' = w^o_{\epsilon_j}F^{-m_j''}_{\epsilon_j}q' +m_j''+ \left [ {\pi\over \epsilon} - n_j + o(1) \right ] $$
As we let $j\to\infty$, the left hand side will converge to $\varphi^\iota p$.  The term  $[\cdots]$ will converge to $-\alpha$.  Thus we conclude that  $w^o_{\epsilon_j}F^{-m_j''}_{\epsilon_j}q' + m_j''$  will converge to $\varphi^\iota p + \alpha$.  By hypothesis, we have $T_\alpha p\in \Sigma_0$, and $\varphi^o$ is a coordinate on $\Sigma_0$.  Thus $\varphi^o$ is a coordinate on $M_\epsilon$ for $\epsilon$ small.  By Proposition 3.7,  $\hat q\mapsto  w^o_{\epsilon_j}F^{-m_j''}_{\epsilon_j}\hat q +m_j''$ gives a uniform family of coordinates on $M_{\epsilon_j}$,  so we conclude that $q'$ must converge to a point $q_0\in\Sigma$.  By the condition that $\varphi^\iota p=\varphi^o q_0 - \alpha$, we conclude that $q_0=T_\alpha p$.
\qed

\begin{theorem}
If $\{\epsilon_j\}$ is an $\alpha$-sequence, then $F^{n_j}_{\epsilon_j}$ converges to $T_\alpha$ uniformly on compact subsets of ${\cal B}$.
\end{theorem}

\noindent{\it Proof.}  We recall that  ${\cal B}$ and $\Sigma$ are invariant in both forward and backward time.  Further $T_{\alpha+1}=T_\alpha\circ f$.  Thus for an arbitrary point $p\in{\cal B}$ and arbitrary $\alpha\in{\bf C}$ we may map $p$ and add an integer to $\alpha$ so that the projection  $\pi p=x$ satisfies $|x+r|<r$, and $T_\alpha p\in \Sigma_0$.  Finally, if we iterate $p$ forward, it will approach $M_0$.  Thus we may also assume that $p\in M_\epsilon$, and now we are in the hypotheses of Proposition 3.8.  \qed

\medskip

\section{Semi-continuity of Julia sets.}  We say that  $\zeta_0\in\bC$ is a periodic point for $h_\alpha$ if $\zeta_j:=h_\alpha^j(\zeta_0)\in\Omega$ for all $j$, and $h_\alpha^n(\zeta_0)=\zeta_0$.  By the chain rule we have $(h_\alpha^n)'(\zeta_0)=\prod_{j=0}^{n-1} h_\alpha'(\zeta_j)$.  We say that $\zeta_0$ is a repelling periodic point if $|(h_\alpha^n)'(\zeta_0)|>1$.  Let $\cR_\alpha$ denote the set of repelling periodic points of $h_\alpha$, and define $J^*(F,T_\alpha)$ to be the closure in $M$ of $(\varphi^o)^{-1}(\cR_\alpha)$.

\begin{theorem}
Let $\zeta_0$ be a repelling (resp.\ attracting) periodic point of period $\mu$ for $h_\alpha$, and let $p_0=(\varphi^o)^{-1}(\zeta_0)\in\cB\cap \Sigma$ be its image.  Then there exists $j_0$ such that for $j\ge j_0$, there is a point $p_j$ near $p_0$, which has period $\nu_j$  for $F_{\epsilon_j}$, with $\epsilon_j={\pi\over j+\alpha}$ and which is a saddle (resp. sink).  Further,  $\nu_j$ divides $j\mu$, and $\nu_j\to\infty$. 
\end{theorem}

\noindent{\it Proof.}  We will prove the repelling case; the attracting case is similar.   If $\zeta_0$ is a repelling periodic point, any the closure of any small disk $\Delta_0$ containing $\zeta_0$ will be contained in $h_\alpha^\mu\Delta_0$.  Let us write $\Delta_0$ for the image of $\Delta_0$ under $(\varphi^o)^{-1}$, and let us consider a neighborhood of $p_0$ which is essentially a product neighborhood, which we may write as $\Delta_0\times\Delta'$.  It follows that $R_\alpha^\mu$ maps $\Delta_0\times\Delta'$ to a disk in $\Sigma_0$ with the following properties:  $R_\alpha^\mu(\partial \Delta_0\times\overline{\Delta'})\cap \overline{\Delta_0\times\Delta'}=\emptyset$, and $\overline{R^\mu_\alpha(\Delta_0\times\Delta')}\cap\overline{(\Delta_0\times\Delta')}\cap \overline{\Delta_0}\times\partial\Delta' = \emptyset$.  By Theorem 3.9, the sequence $F^j_{\epsilon_j}$ converges uniformly on compacts to $R_\alpha$, and thus $F^{j\mu}_{\epsilon_j}$ converges uniformly on $\overline{\Delta_0\times\Delta'}$ to $R_\alpha^\mu$.  It follows that $F^{j\mu}_{\epsilon_j}$ has the same mapping properties on the product $\Delta_0\times\Delta'$.   Thus $F^{j\mu}_{\epsilon_j}$ has a saddle point $p_j$ in $\Delta_0\times\Delta'$.  The period $\nu_j$ of $p_j$ must divide $j\mu$.  Since $\Delta_0\times\Delta'$ can be taken arbitrarily small, we see that the $p_j$ will converge to $p_0$.  Finally, $\nu_j$ cannot have a bounded subsequence, or else $p_0$ would be periodic for $F_0$.  But this is impossible since $p_0\in\cB$.
\qed

For an automorphism $F$, we define $J^*=J^*(F)$ to be  the closure of the saddle periodic points of $F$.  In general the sets $J^*(F_\epsilon)$ are lower semicontinuous as a function of $\epsilon$.   Since $J^*(F)\cap \cB=\emptyset$, the following result gives a lower estimate for the discontinuity of the sets $J^*$ and gives a proof of Theorem 1:

\begin{theorem}
If $\{\epsilon_j\}$ is an $\alpha$-sequence, then $\liminf_{j\to\infty}J^*(F_{\epsilon_j})\supset J^*(F,T_\alpha)$.
\end{theorem}

\noindent{\it Proof.}  Let $p_0$ be a periodic point in $\cJ^*_\alpha$.  It will suffice to show that for every $\epsilon>0$ there is a $j_0$ such that for $j\ge j_0$ there is a saddle point $p_j$ for $F_{\epsilon_j}$ which is within $\epsilon$ of $p_0$.  This property is given by the previous Theorem.   \qed

Now let us suppose that $M=\bC^2$ and  $F:M \to M$ is a composition $F = F_1\circ\cdots\circ F_k$, where $F_j(x,y) = (y, P_j(y)-\delta_j x)$ is a  generalized H\'enon map.  We define $K^\pm$ as the points with bounded forward/backward orbits.  
The set $K:=K^+\cap K^-$ is bounded and contains $\cB\cap\Sigma$, so we have:  
\begin{proposition}
If $F$ is a  composition of generalized H\'enon maps, then $\Sigma\not\subset\cB$.  Furthermore, every component of $\cB\cap\Sigma$ is conformally equivalent to the disk.
\end{proposition}

\noindent{\it Proof.}  We have $\cB\subset K^+$ and $\Sigma\subset K^-$, so $\cB\cap\Sigma\subset K$, which is bounded.  Thus $\cB\cap\Sigma$ cannot be uniformized by $\bC$.   By Proposition 1.4 it is simply connected, so it must be a disk.  \qed

\bigskip
\centerline{\includegraphics[width=3in]{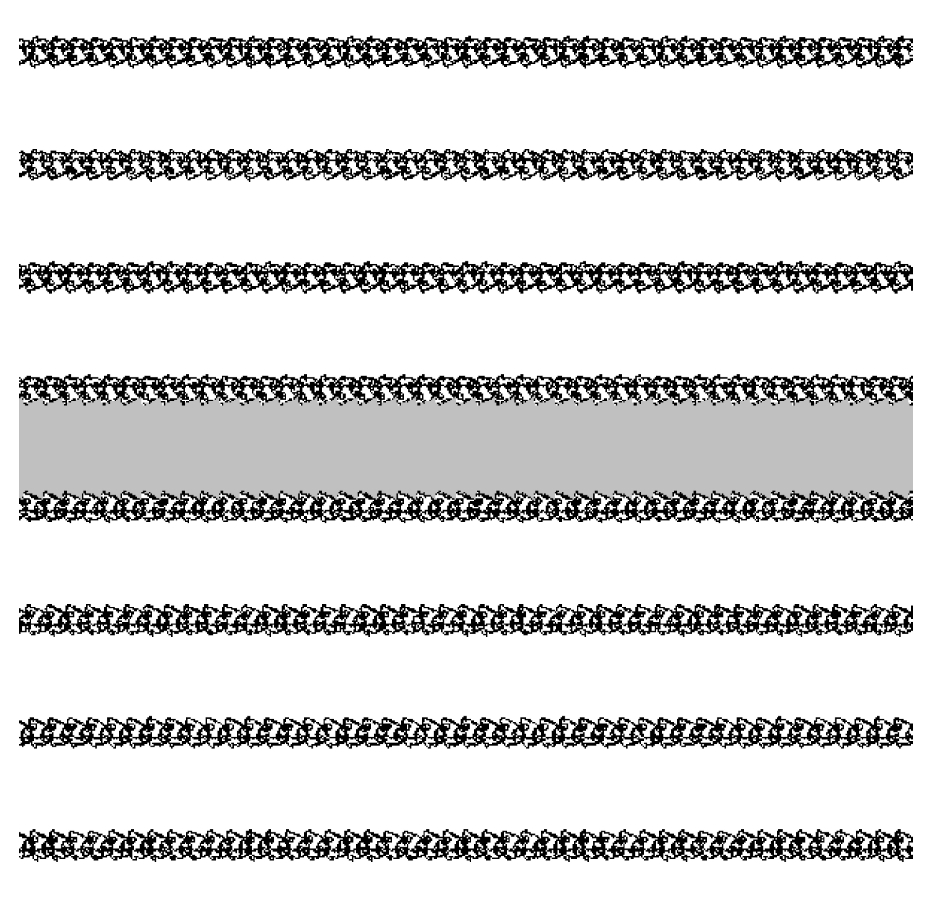} \ \ \  \includegraphics[width=3in]{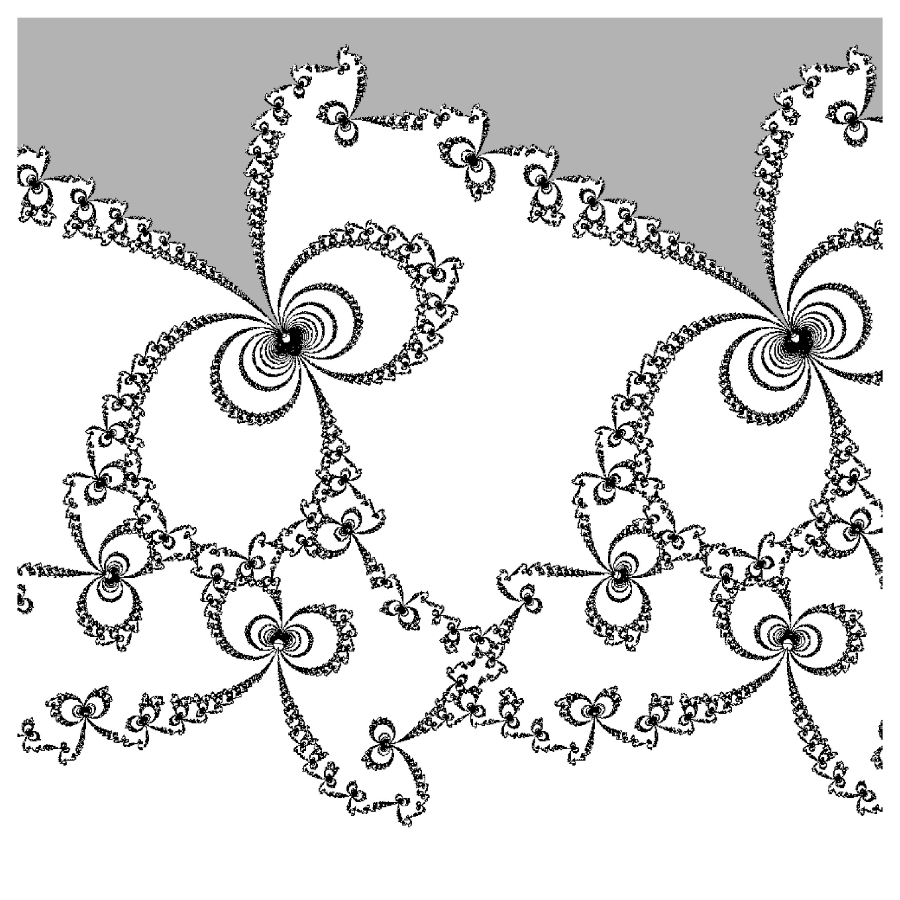}   }


\centerline{Figure 5.  Slice $K^+(F,T_\alpha)\cap \Sigma$, $a=.3$, $\alpha=0$: 44 periods, $|\Im\zeta|<22$ (left); detail (right). }
\bigskip

We define an analogue of a ``Julia-Lavaurs set'':  we define $K^+(F,T_\alpha)$ to be the set of points $p\in K^+(F)$ which satisfy one of the following two properties:
\begin{enumerate}
\item[(i)]  $T^n_\alpha(p)$ is defined and belongs to $\cB$ for all $n\ge0$.
\item[(ii)]  There is an integer $n\ge0$ such that $T^k_\alpha(p)\in\cB$ for $k\le n-1$, and $T^n_\alpha(p)\in K^+-\cB$.
\end{enumerate}

Thus the complement of $K^+(F,T_\alpha)$ consists of the points satisfying the condition: there is an $n\ge0$ such that $T^k_\alpha(p)\in\cB$ for $k\le n-1$ and that $T^n_\alpha(p)\notin K^+$.  It is immediate from the definition that
$$K^+(F,T_\alpha)-\cB= K^+(F)-\cB\subset K^+(F,T_\alpha)\subset K^+(F).$$

\begin{theorem}
There exist $p\in \cB$ and $\alpha'$ such that  $p\in J^*(F_0,T_\alpha)$.    Further, for each $p\in\cB$ there is an $\alpha'$ such that $p\notin K^+(F_0,T_{\alpha'})$.     In particular, we may choose $p\in \cB$, $\alpha$ and  $\alpha'$ such that $p\in J^*(F_0,T_\alpha)$, but $p\notin K^+(F_0,T_{\alpha'})$.
\end{theorem}

\noindent{\it Proof.}  Since $F_0$ is a H\'enon map, there is a point $q\in\Sigma -K^+$.  Thus, given $p$, we choose $\alpha$ so that $\alpha = \varphi^o(q) - \varphi^\iota(p)\in{\bf C}$.  It follows that $p\notin K^+(F_0,T_\alpha)$. 

Next, we consider the partially defined map $h_0:=\varphi^\iota\circ(\varphi^o)^{-1}:{\bf C}\to{\bf C}$.  By \S1, we know that $h_0$ is a well-defined map of the cylinder ${\bf C}/{\bf Z}$, and  $h_0'(\zeta)\to1$ as $\zeta$ approaches either end of the cylinder.  By Proposition 4.3, $h_0$ cannot be holomorphic in a neighborhood of either end of the cylinder.  Thus we must have both $|h'_0|>1$ and $|h'_0|<1$ at points near either end of the cylinder.  Chose a point $\zeta_0$ such that $|h'_0(\zeta_0)|>1$.  Then we may choose $\alpha\in{\bf C}$ such that $h_\alpha(\zeta_0)=\tau_\alpha(h_0(\zeta_0)) = \zeta_0$.  It follows that  $\zeta_0$ is a repelling fixed point for $h_\alpha$.  Thus $(\varphi^o)^{-1}(\zeta_0)\in J^*(F_0,T_\alpha)$.   \qed

\centerline{\hbox{ \includegraphics[width=3in]{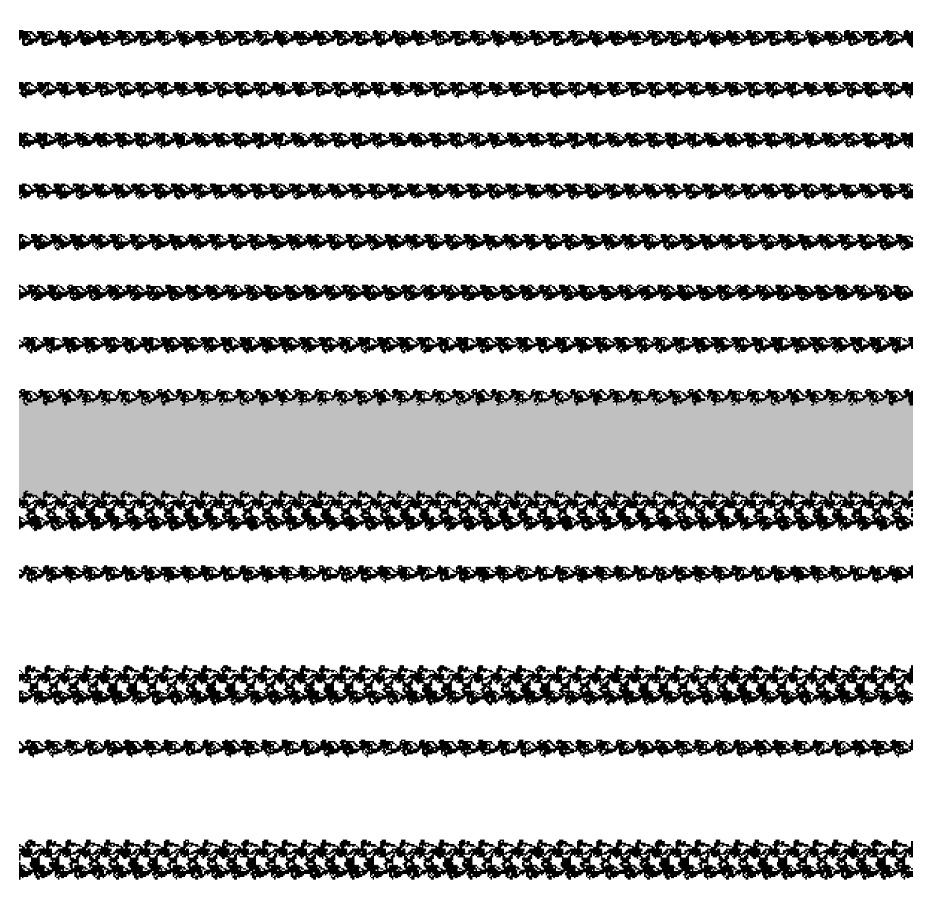}  
	\includegraphics[width=3in]{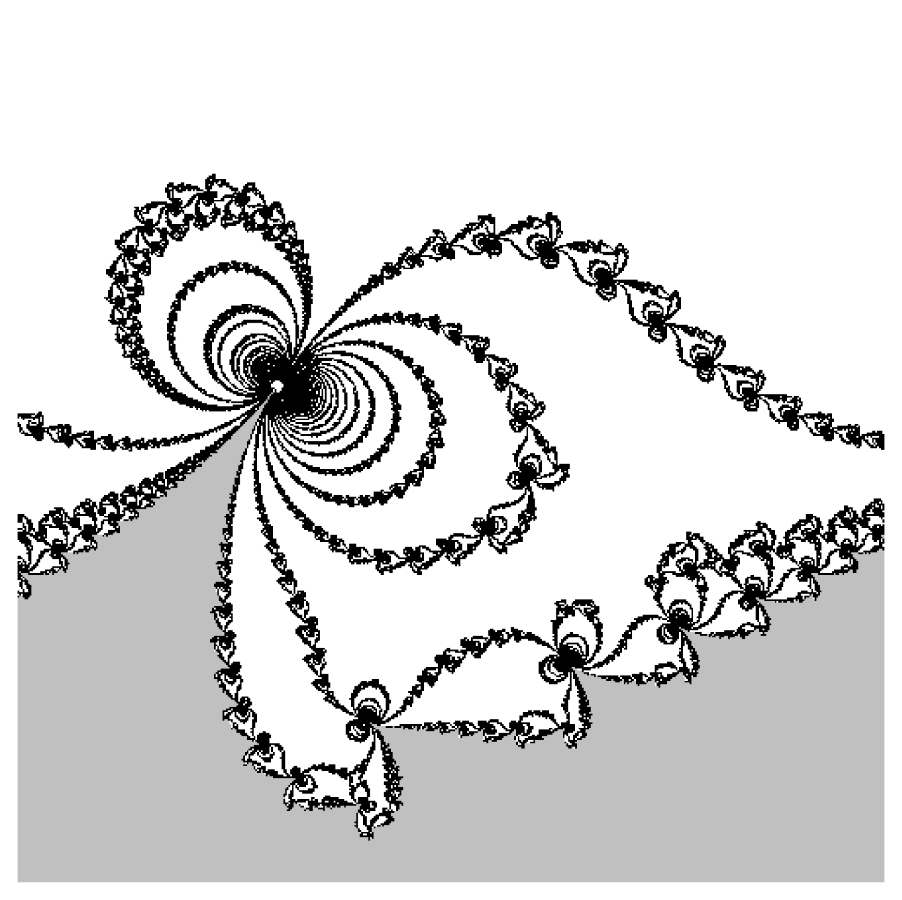} } }

%
%
\medskip

\centerline{Figure 6a.  Slices $K^+(F,T_\alpha)\cap \Sigma$; $a=.3$, $\alpha=\pi i$; 44 periods, $|\Im\zeta|<22$ (left), detail  (right). }

\bigskip

To illustrate $K^+(F,T_\alpha)$ graphically, we return to the H\'enon family defined in (1.1).  The pictures in Figure 5 correspond to Figure 2.  That is, they are slices of $K^+(F,T_\alpha)\cap\Sigma$, with the values $a=.3$ and $\alpha=0$, which corresponds to real $\epsilon$.  The gray region is the complement of $K^+(F)$, the set $K^+(F,T_\alpha)$ is black, and $K^+(F)-K^+(F,T_\alpha)$ is white.  All pictures are invariant under the translation $\zeta\mapsto\zeta+1$. The viewboxes on the left hand sides of Figures  5 and 6 are taken to be symmetric around the real axis $\{\Im\zeta=0\}$;  the viewboxes are taken to have side $=44$ in order to show what happens when $\Im\zeta$ is large.  We see a number of horizontal ``chains'' in the left hand pictures in Figures 5 and 6.  In the upper half of each of these pictures, the map $h_\alpha$ acts approximately as a vertical translation, moving each chain to the one below it, until it reaches the chain just above the gray region, which corresponds to the complement of ${\cal B}$.  By (1.7) the amount of vertical translation in the upper region is approximately  $c_0^+\approx -5.83$.  This fits with the height of the box in Figure 5, since there are $8\approx 44/5.83$ horizontal strips.  In Figure 6, the vertical translation in the upper part is $c_0^+ + \Im\alpha\approx -2.69$.  In the chains bordering the complement of the basin, the map is not like a translation and is more complicated.  The bottom half of the left hand side of Figure 5 and 6 is analogous, with the approximate translation near the bottom of the figures being approximately $c_0^-+\Im\alpha$.  In fact, the symmetry in Figure~5 comes because $h_\alpha$ commutes with complex conjugation.  The pictures on the right of Figures 5 and 6a,b give a detail from the edge of the gray region, spanning a little more than 1 period.  The implosion phenomenon corresponding to Figure~6 is given in Figures 7b,c, where we see all three pictures from Figure 6.

%

%


%

%

\centerline{ \includegraphics[width=2.2in]{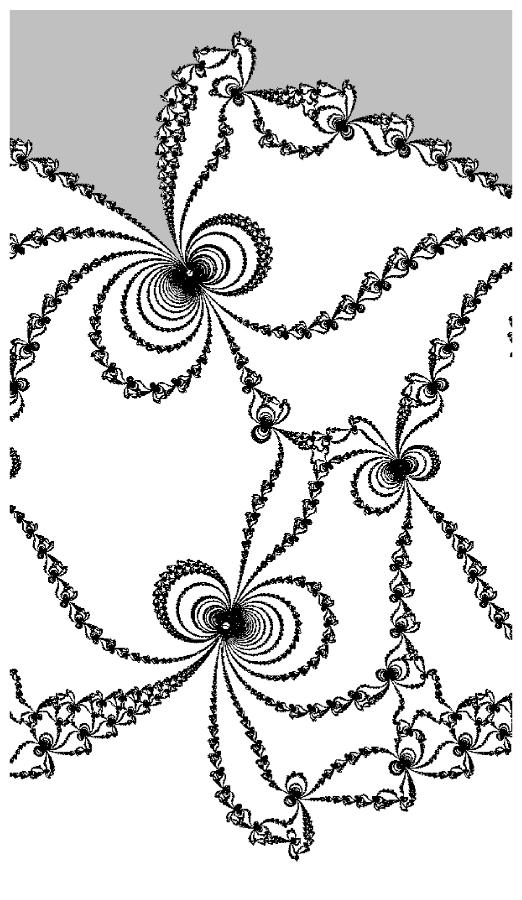}    }

\centerline{Figure 6b. Detail of 6a: one period from top of the bottom component of $\Sigma\cap{\cal B}$ }
\medskip

\begin{proposition}
$K^+(F,T_\alpha)\cap{\cal B}$ is a union of fibers $\{\varphi^\iota=const\}$, and $K^+(F,T_\alpha)\cap{\cal B}\ne{\cal B}$.
\end{proposition}

\noindent{\it Proof.}  The first statement follows from the definition $T_\alpha:=(\varphi^o)^{-1}\circ\tau_\alpha\circ\varphi^\iota$.  For the second statement, we recall that $\Sigma\not\subset K^+$.  So choose a point $p_0=(\varphi^o)^{-1}(\zeta_0)\notin K^+$.  It follows that the fiber $\{\varphi^\iota=\zeta_0-\alpha\}$ is mapped to $p_0$.  Thus this fiber is outside of $K^+(F,T_\alpha)$.   \qed

\noindent{\bf Remark. } Let $p\in \Sigma\cap K^+(F,T_\alpha)\cap \cB$ be a point which is not critical for $h_\alpha$ (which means that $\Sigma$ is not tangent to the fibers of $\varphi^\iota$ at $p$).  Then in a neighborhood of $p$,  $K^+(F,T_\alpha)$ will be a product of a disk with the slice $K^+(F,T_\alpha)\cap \Sigma$.

We also have the following elementary observation:

\begin{proposition}
$K^+(F,T_\alpha)=FK^+(F,T_\alpha)=K^+(F,T_{\alpha+1})$, so $K^+(F,T_\alpha)$ depends only on the equivalence class of $\alpha$ modulo $\bZ$.
\end{proposition}

\begin{proposition}
In general, $F\mapsto K^+(F)$ is upper semicontinuous: $\limsup_{\epsilon_j\to0}K^+(F_{\epsilon_j})\subset K^+(F_0)$, and $F\mapsto J^+(F)$, $F\mapsto J(F)$ are lower semicontinuous, $\liminf_{j\to\infty} J^+(F_\epsilon)\subset J^+(F_0)$.  Similar statements hold for $K^-$ and $J^-$.
\end{proposition}

\noindent{\it Proof.  }  This follows because the Green function $G^+\epsilon$ is continuous and depends continuously on $\epsilon$. Thus $K^+_\epsilon=\{G^+_\epsilon=0\}$ is upper semicontinuous.  On the other hand the measure $\mu_\epsilon$ and the current  $\mu^+_\epsilon$ depend continuously  on $\epsilon$, thus their supports $J^+={\rm supp}(\mu^+_\epsilon)$ and $J={\rm supp}(\mu_\epsilon)$ are lower semicontinuous. \qed
 
The following gives a sharpening of the semicontinuity and gives a proof of Theorem~2:
\begin{theorem}
If $\{\epsilon_j\}$ is an $\alpha$-sequence, then ${\cal B}\cap \limsup_{j\to\infty} K^+(F_{\epsilon_j})\subset K^+(F,T_\alpha)$.
\end{theorem}

\noindent{\it Proof of Theorem 2. }  Let us choose a point $p\notin K^+(F,T_\alpha)$.  Thus there exists an $m$ such that $T^m_\alpha(p)\notin K^+(F_0)$.  It will suffice to show that $p\notin K^+(F_{\epsilon_j})$ for large $j$.  By Theorem 3.9, it follows that $f^{n_j m}_{\epsilon_j}p$ is approximately $T^m_\alpha (p)$, and thus  $f^{n_j m}_{\epsilon_j}p\notin K^+(F_0)$.  By the semicontinuity of $K^+$, it follows that $f^{n_j m}_{\epsilon_j}p\notin K^+(F_{\epsilon_j})$.  Thus $p\notin K^+(F_{\epsilon_j})$.   \qed



Since $F$ is a polynomial automorphism, the Jacobian is constant, and since the parabolic point is semi-attracting, $F$ contracts volume.  A consequence is $J^-(F)=K^-(F)$.  We define 
\begin{align*}
K(F,T_\alpha):&=J^-(F)\cap K^+(F,T_\alpha)=K^-(F)\cap K^+(F,T_\alpha),\cr
J(F,T_\alpha):   &= J^-(F)\cap \partial K^+(F,T_\alpha).
\end{align*}
Thus we have $\cJ^*_\alpha\subset J(F,T_\alpha)$.  Since $K^-(F_\epsilon)$ is upper semicontinuous, we have:

\begin{corollary}
If $\{\epsilon_j\}$ is an $\alpha$-sequence, then ${\cal B}\cap \limsup_{j\to\infty} K(F_{\epsilon_j})\subset K(F,T_\alpha)$.
\end{corollary}

\vbox{
\bigskip

\centerline{  \includegraphics[height=2.5in]{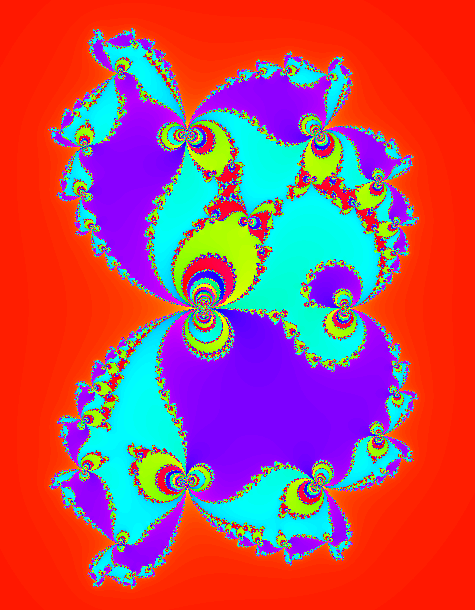} \ \ \ \ \  \ \ \includegraphics[height=2.5in]{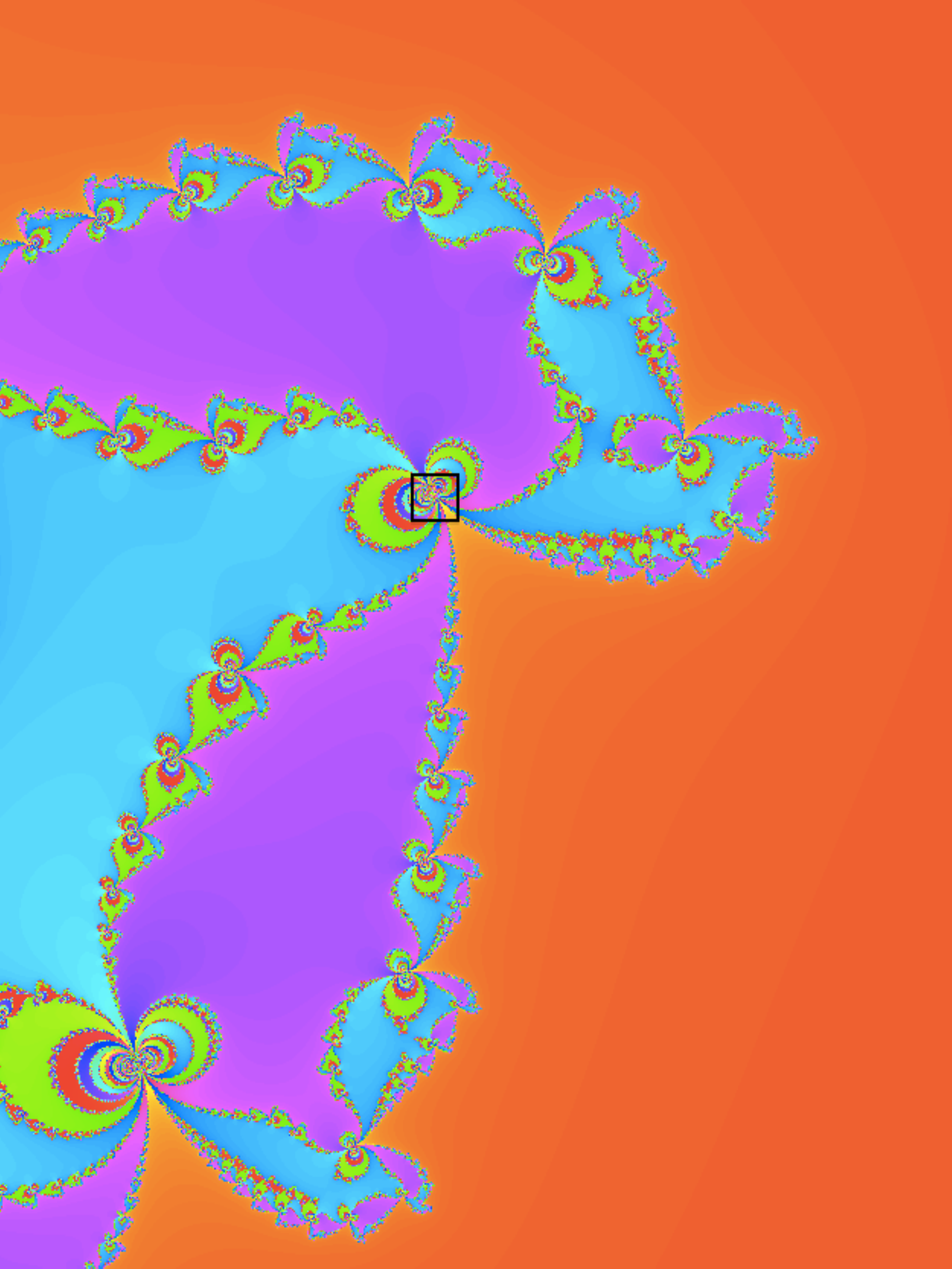} }
\medskip

\centerline{Figure 7a. Slices of $K^+$ for $f$ in (1.1), with $a=.3$, $\epsilon=\pi/(n- i \alpha)$, $n=1000$, $\alpha=4.3$:}  
\vskip0in

\centerline{Linear slice $K^+\cap T$ (left),  unstable slice $K^+\cap W^u(q)$ (right)} }

 \bigskip

\centerline{\ \  \includegraphics[height=2.5in]{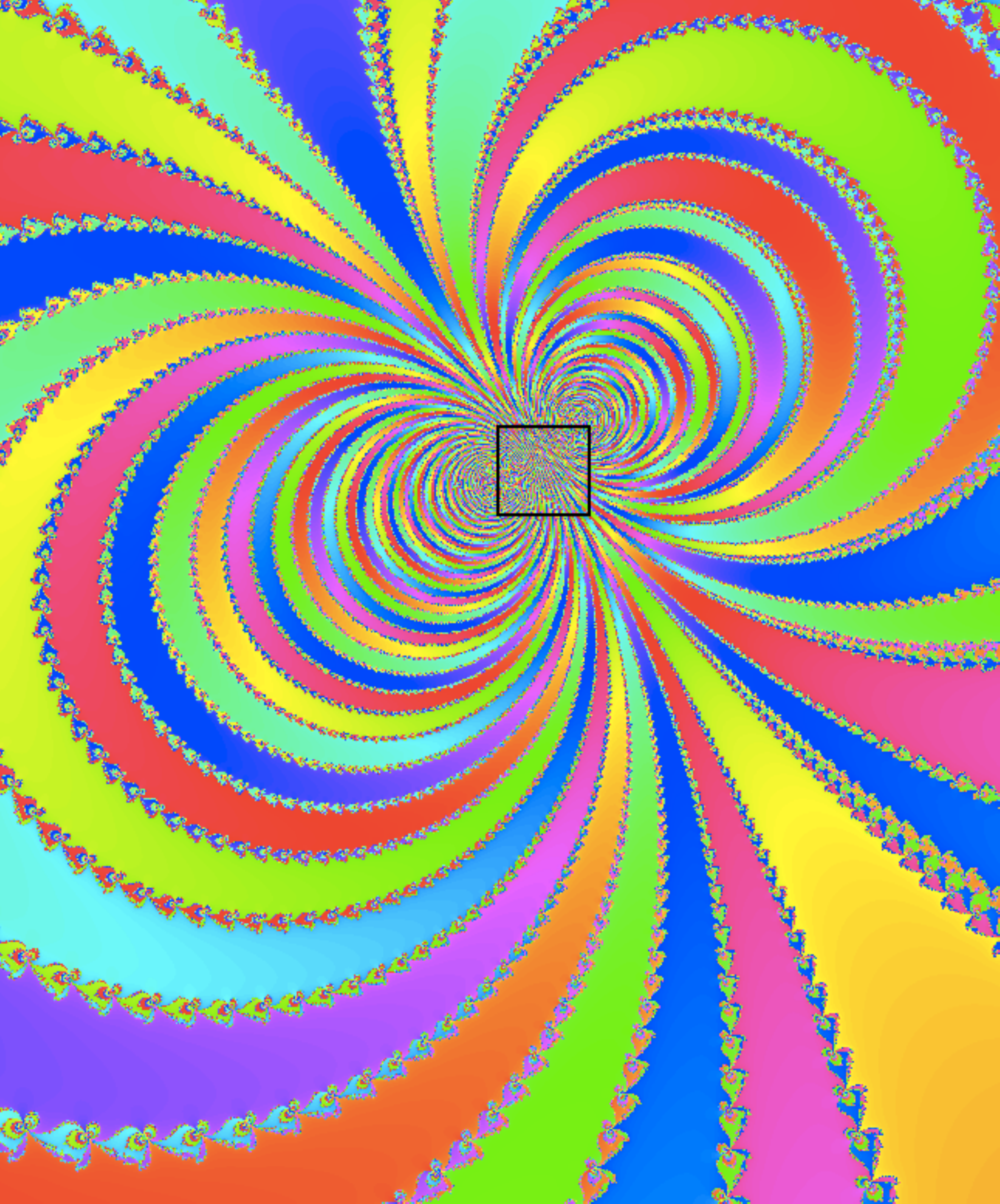} \ \ \ \ \ \  \includegraphics[height=2.5in]{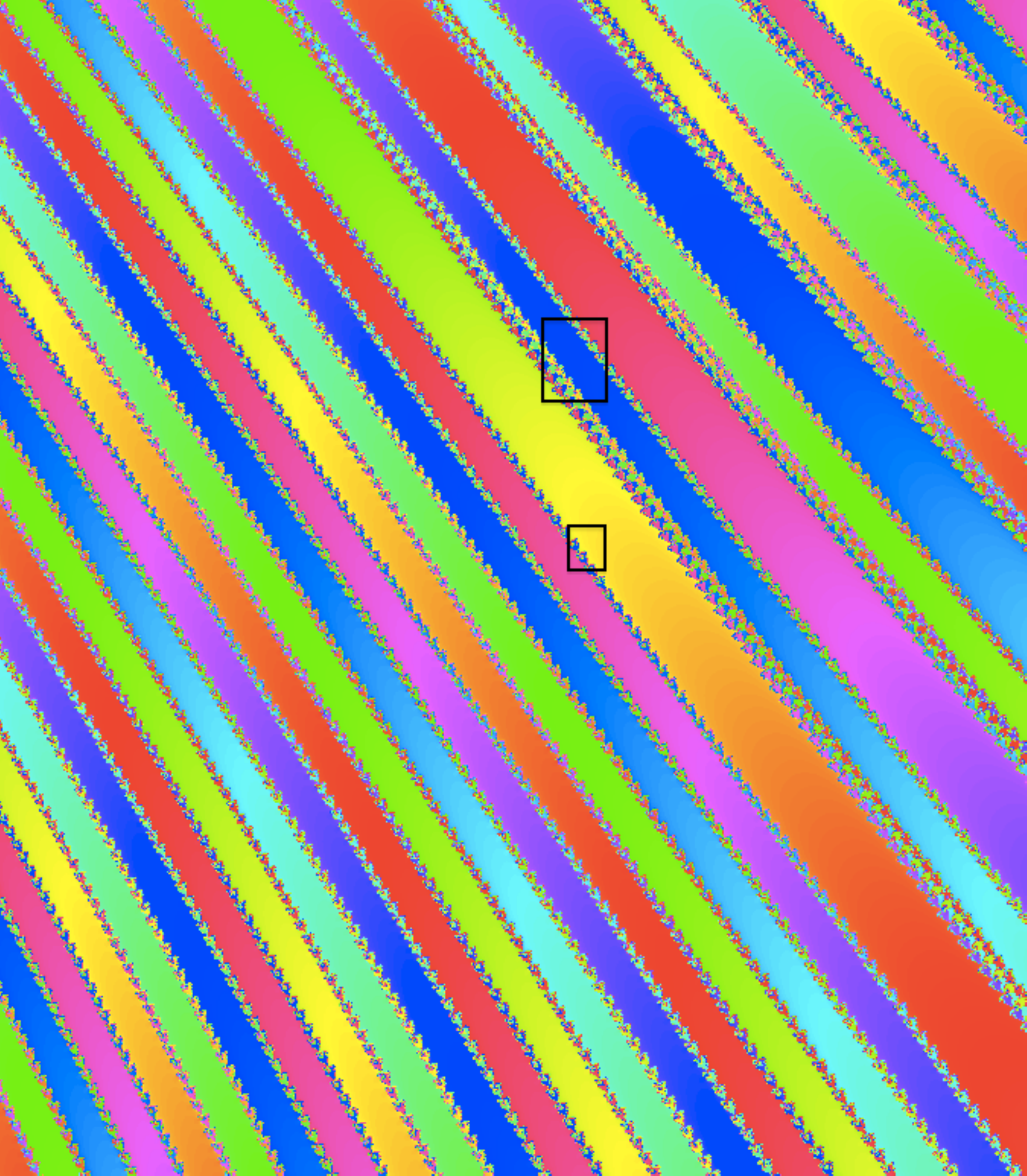} }

\centerline{Figure 7b. Slices of $K^+$ for $f$ in (1.1), with $a=.3$, $\epsilon=\pi/(n- i \alpha)$, $n=1000$, $\alpha=4.3$:  Further zooms.}  

\bigskip 

\centerline{ \includegraphics[height=2.7in]{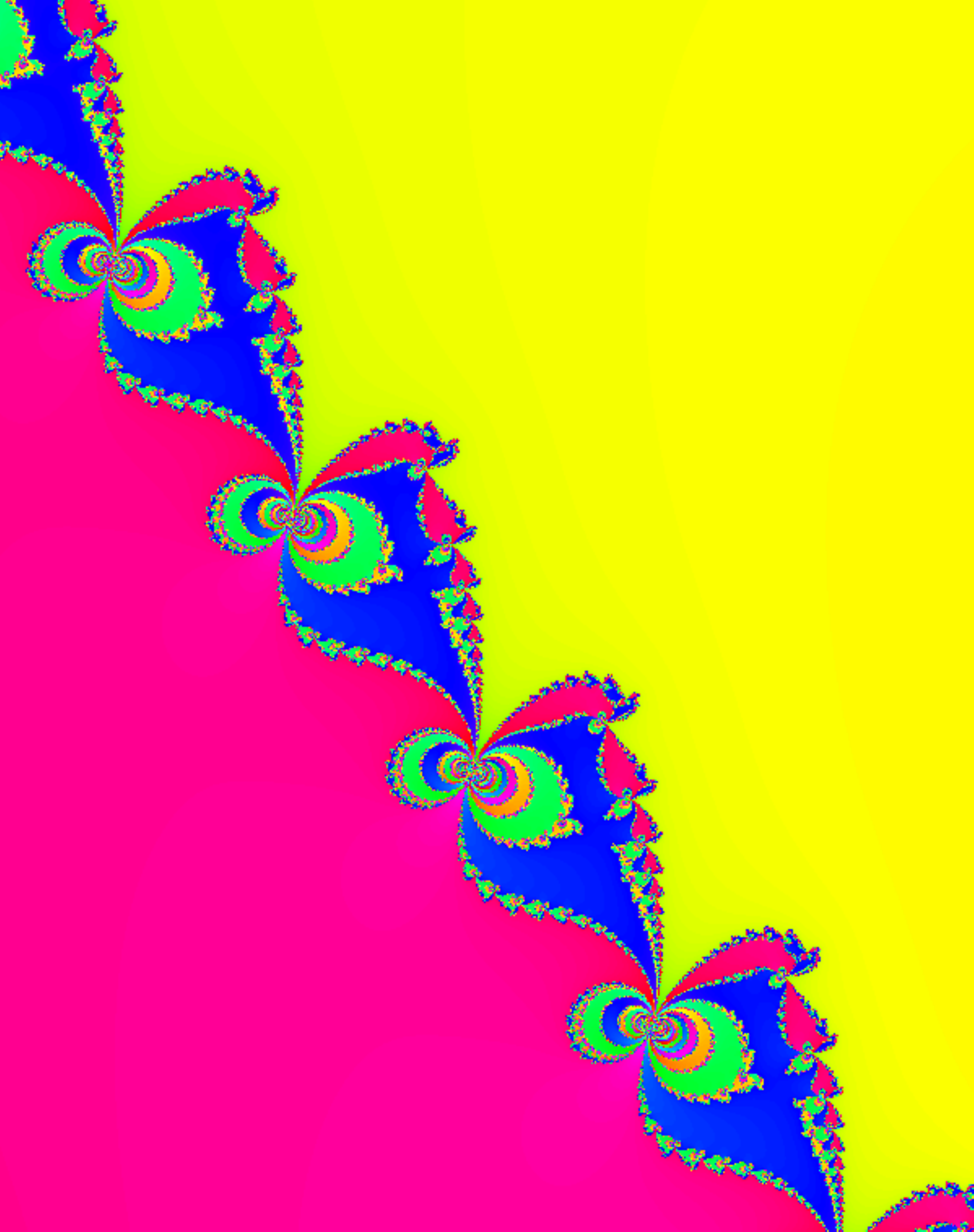} \ \ \ \ \ \ \              \includegraphics[height=2.7in]{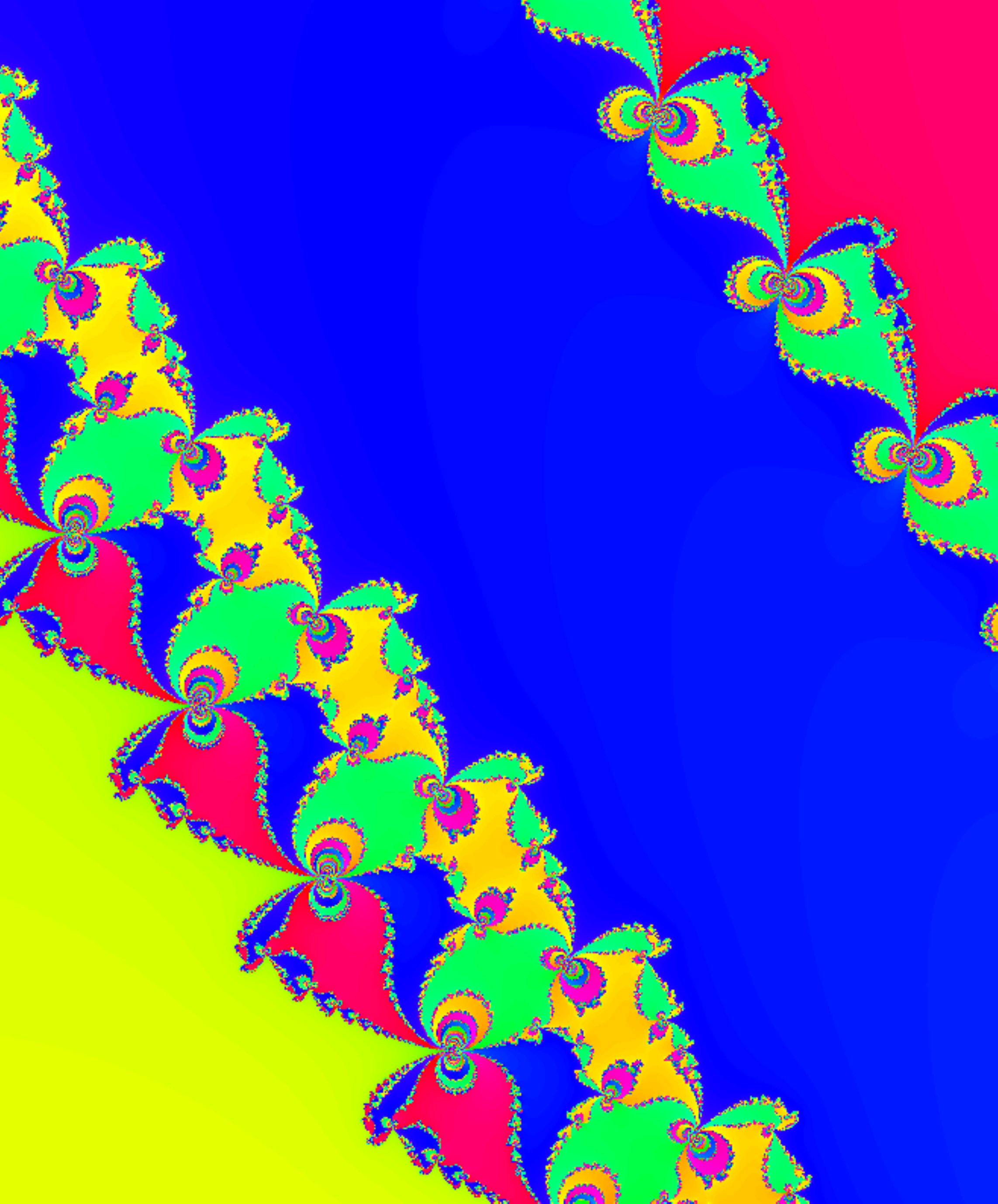}    }

\centerline{Figure 7c. Details of right side of Figure 7b.}  

\bigskip

\bigskip

\centerline{\bf References}
\medskip

\begin{enumerate}
\item[{[BS1]}]  E. Bedford and J. Smillie, Polynomial diffeomorphisms of ${\bf C}^2$: currents, equilibrium measure and hyperbolicity. Invent. Math. 103 (1991), no. 1, 69--99.

\item[{[BS7]}]  E. Bedford and J. Smillie,  Polynomial diffeomorphisms of $\bC^2$. VII. Hyperbolicity and external rays. Ann. Sci. \'Ecole Norm. Sup. (4) 32 (1999), no. 4, 455--497.



\item[{[D]}]  A. Douady, Does a Julia set depend continuously on the polynomial? {\sl Complex dynamical systems (Cincinnati, OH, 1994)}, 91--138, Proc.\ Sympos.\ Appl.\ Math., 49, Amer.\ Math.\ Soc., Providence, RI, 1994.

\item[{[DH]}] A. Douady and J. Hubbard, \'Etude dynamique des polyn\^omes complexes. Partie I.  Publications Math\'ematiques d'Orsay, 84-2. Universit\'e de Paris-Sud, 1984.

\item[{[DSZ]}]  A. Douady,  P. Sentenac, and M.  Zinsmeister,   Implosion parabolique et dimension de Hausdorff. C. R. Acad.\ Sci.\ Paris S\'er.\ I Math. 325 (1997), no. 7, 765--772. 

\item[{[FS]}]  J.-E. Forn\ae ss and N. Sibony, Complex H\'enon mappings in ${\bf C}^2$ and Fatou-Bieberbach domains. Duke Math.\  J. 65 (1992), no.\  2, 345--380.

\item[{[H]}] M. Hakim,  Attracting domains for semi-attractive transformations of $\bC^p$. Publ. Mat. 38 (1994), no. 2, 479--499.

\item[{[HO]}]  J.H. Hubbard and R. Oberste-Vorth, H\'enon mappings in the complex domain. I. The global topology of dynamical space. Inst.\ Hautes \'Etudes Sci.\ Publ.\ Math.\ No. 79 (1994), 5--46.

\item[{[HPS]}]  M. Hirsch, C. Pugh, and M. Shub,  Invariant manifolds. {\sl Lecture Notes in Mathematics}, Vol.\ 583. Springer-Verlag, Berlin-New York, 1977. ii+149 pp. 

\item[{[L]}]  P. Lavaurs,  Syst\`emes dynamiques holomorphiques:  explosion de points p\'eriodiques.  Th\`ese, Universit\'e Paris-Sud, 1989.

\item[{[Mc]}] C. McMullen,  Hausdorff dimension and conformal dynamics. II. Geometrically finite rational maps. Comment.\ Math.\ Helv.\ 75 (2000), no. 4, 535--593.

\item[{[Mi]}]  J. Milnor,  {\sl Dynamics in One Complex Variable}, Annals of Math.\ Studies, 2006.

\item[{[O1]}]  R. Oudkerk,  The parabolic implosion: Lavaurs maps and strong convergence for rational maps. {\sl Value distribution theory and complex dynamics (Hong Kong, 2000)}, 79--105,  Contemp. Math., 303, Amer.\ Math.\ Soc., Providence, RI, 2002.

\item[{[O2]}]  R. Oudkerk, Parabolic implosion for $f_0(z) = z + z^{\nu+1}+{\cal O}(z^{\nu+2})$, Thesis, U. of Warwick, 1999.

\item[{[S]}]  M. Shishikura,   Bifurcation of parabolic fixed points. {\sl The Mandelbrot set, theme and variations}, 325--363, London Math.\ Soc.\ Lecture Note Ser., 274, Cambridge Univ.\ Press, Cambridge, 2000.

\item[{[U1]}]  T. Ueda,  Local structure of analytic transformations of two complex variables. I. J. Math. Kyoto Univ. 26 (1986), no. 2, 233--261.

\item[{[U2]}]  T. Ueda,  Local structure of analytic transformations of two complex variables. II. J. Math. Kyoto Univ. 31 (1991), no. 3, 695--711.





\item[{[Z]}]  M. Zinsmeister,  Parabolic implosion in five days, Notes from a course given at Jyvaskyla, September, 1997.

\end{enumerate}

\bigskip
\rightline{Eric Bedford}

\rightline{Indiana University}

\rightline{Bloomington, IN 47405}

\rightline{\tt bedford@indiana.edu}

\bigskip
\rightline{John Smillie}

\rightline{Cornell University}

\rightline{Ithaca, NY 14853}

\rightline{\tt smillie@math.cornell.edu}

\bigskip
\rightline{Tetsuo Ueda}

\rightline{Kyoto University}

\rightline{Kyoto 606-8502, Japan}

\rightline{\tt ueda@math.kyoto-u.ac.jp}

\vfill\eject

\end{document}